\documentclass{article}
\usepackage[utf8]{inputenc}

\usepackage{latexsym}
\usepackage{amssymb}
\usepackage{epsfig}
\usepackage{psfig}
\usepackage{amssymb}
\usepackage{amsmath}
\usepackage{amsthm}
\usepackage{graphicx}
\usepackage{color}
\usepackage{xspace}
\usepackage{amssymb}
\usepackage{amsfonts}
\usepackage{latexsym}
\usepackage{epsfig}
\usepackage{psfig}
\usepackage{xspace}
\usepackage{color}
\usepackage{amstext}
\usepackage{rotating}
\usepackage{amsmath}
\usepackage{amsthm}

\renewcommand{\qedsymbol}{$\dashv$}

\newcommand{\sucs}{\prec}

\newcommand{\mcs}{\ensuremath{\textup{MCS}}\xspace}

\newcommand{\eqbydef}{:=}

\newcommand{\extil}[1]{\ensuremath{\textup{\textbf{IL}}{\sf\ensuremath{#1}}}\xspace}

\newcommand{\il}{{\ensuremath{\textup{\textbf{IL}}}}\xspace}
\newcommand{\gl}{{\ensuremath{\textup{\textbf{GL}}}}\xspace}

\newcommand{\ilw}{\extil{W}}
\newcommand{\ilp}{\extil{P}}

\newcommand{\ilx}{\extil{X}}

\theoremstyle{plain}
\newtheorem{theorem}{Theorem}[section]
\newtheorem{lemma}[theorem]{Lemma}

\newtheorem{notation}[theorem]{Notation}

\theoremstyle{definition}
\newtheorem{definition}[theorem]{Definition}

\newtheorem{remark}[theorem]{Remark}

\theoremstyle{remark}

\newcommand{\sassuring}[1]{{\ensuremath{\prec_{#1}}}\xspace}

\newcommand{\principle}[1]{{\ensuremath{\sf #1}}\xspace}
\newcommand{\pair}[1]{{\ensuremath{\langle #1 \rangle}}\xspace}
\newcommand{\boxset}[2]{{\ensuremath{{#2}^{\Box}_{#1}}}\xspace}

\usepackage{amsmath}
\usepackage{tikz}
\usepackage{mathdots}
\usepackage{yhmath}
\usepackage{cancel}
\usepackage{color}
\usepackage{siunitx}
\usepackage{array}
\usepackage{multirow}
\usepackage{amssymb}
\usepackage{gensymb}
\usepackage{tabularx}
\usepackage{booktabs}
\usepackage{wasysym}
\usetikzlibrary{fadings}

\newcommand{\illuka}[1]{\text{\textbf{IL}\textsf{#1}}}

\newcommand{\kgen}[1]{\text{(\textsf{{#1}})\textsubscript{gen}}}

\usepackage{enumitem}

\author{
  Evan Goris\footnote{Independent scholar}\\
  Marta Bílková\footnote{Institute of Computer Science of the Czech Academy of Sciences.}\\
  Joost J. Joosten\footnote{Department of Philosophy, University of Barcelona. Supported by the Spanish Ministry of Science and Universities under grant number RTC-2017-6740-7, Spanish Ministry of Economy and Competitiveness under grant number FFI2015-70707P and the Generallitat de Catalunya under grant number 2017 SGR 270.}\\
  Luka Mikec\footnote{Department of Mathematics, Faculty of Science, University of Zagreb. 
  Supported by Croatian Science Foundation (HRZZ) under the projects UIP-05-2017-9219 and IP-01-2018-7459.}
}

\title{Assuring and critical labels for relations between maximal consistent sets for interpretability logics}

\begin{document}
\maketitle

\begin{abstract}
The notion of a \emph{critical successor} \cite{jovelt90} in relational semantics has been central to 
    most classic modal completeness proofs in interpretability logics. 
In this paper we shall work with a more general notion, that of an \emph{assuring successor}.
As we shall see, this will enable more concisely formulated completeness proofs,
    both with respect to ordinary and generalised Veltman semantics.
Due to their interesting theoretical properties, we will devote some space
    to the study of a particular kind of assuring labels,
    the so-called \emph{full labels} and \emph{maximal labels}.
After a general treatment of assuringness, we shall apply it to obtain 
    certain completeness results.
Namely, we give another proof of completeness of 
    \extil{W} w.r.t.\ ordinary semantics
    and of \extil{P} w.r.t.\ generalised semantics.

\end{abstract}

\section{Introduction}
This paper is about a technical aspect of interpretability logics. Interpretability logics are propositional modal logics that naturally extend provability logic. 

The provability logic \gl is a propositional modal logic with a unary modality $\Box$ that describes the notion of formal provability. Consequently, the dual modality $\Diamond$ refers to consistency. 
The logic \gl comes with a natural Kripke semantics where the truth conditions concerning the modalities are modeled using a binary accessibility relation usually denoted by $R$.  

It is well-known that \gl is complete with respect to its relational semantics (\cite{segerberg1971}). The modal completeness proof is as usual performed by building a sort of term model. That is to say, we build semantics out of syntax. Thus, one considers maximal \gl-consistent sets which will be worlds in the Kripke model. The $R$ relation between maximal consistent sets is defined in such a way that the resulting structure yields a \gl-model. Here and elsewhere, by $S$ being a \textit{maximal consistent set} we understand that $S \nvdash \bot$ and $S' \vdash \bot$ for any proper superset $S'$ of $S$. Here $\vdash$ depends on the logic \textbf{L} which will always be clear from the context (in this paragraph, it is the logic \gl). We write $S \vdash A$ if there is a finite sequence of formulas ending with $A$, each of which is a theorem of \textbf{L}, or contained in $S$, or a consequence by modus ponens of some preceding formulas.
From now on we will write MCS as short-hand for \emph{maximal consistent set}.

Provability logic describes, in a precise sense \cite{solo76}, all the behavior about formal provability of a theory that can be proven by that particular theory. In a similar fashion, interpretability logics describe the provable behaviour of relativised interpretability. Now, a binary modality $\rhd$ is used where the intended reading of $\varphi \rhd \psi$ is that some base theory $T$ together with (the arithmetical reading of) $\varphi$ interprets $T$ together with (the arithmetical reading of) $\psi$.\footnote{That $U$ interprets $V$ means that there is an interpretation of $V$ in $U$, i.e.~there is a translation from the language of $V$ to the language of $U$ which preserves both structure and provability. There are different notions of interpretability, many of them compatible with interpretability logics. We will not need further details for this paper as here we are only concerned with modal aspects of interpretability logics. For a precise definition of an interpretation see e.g.~\cite{vis97}.} By doing so, we will see that $\Box \varphi$ is equivalent to $\neg \varphi \rhd \bot$ so that interpretability logics indeed naturally extend provability logic.

Whereas the logic of provability is very stable and the same for basically any sound and strong enough theory (see e.g.~\cite{JonghJumeletMontagna:1991}), the situation with interpretability logics differs a lot. Different theories have different interpretability logics (see \cite{Visser:1990:InterpretabilityLogic, Zambella:1992:Interpretability, shav88, bera90, BilkovaJonghJoosten:2009:PRA,  JoostenIcard:2012:RestrictedSubstitutions}) which make them interesting to study.

Interpretability logics also come with a Kripke-like semantics \cite{jovelt90}. Their models are called \emph{Veltman models} and they naturally extend models for provability. The truth conditions for the binary modality $\rhd$ is now governed by a ternary relation $S$ between worlds in a Veltman model. 

For various interpretability logics, we also know completeness w.r.t.\ the respective class of Veltman models. Again, completeness proofs proceed by constructing a sort of term model.

Since interpretability logics extend provability logic we will again have a binary accessibility relation between MCSs. 
However, it turns out to be much more difficult to define the ternary 
relation $S$. The reason is that a single MCS may be needed in various roles now to validate different formulas (we shall see a concrete example later in Figure \ref{fig:no_maximum_ex} and clarify this slightly below). Thus, there may be a need to have various copies of a single MCS that occur in different parts of the model to fulfill the different roles of the MCS.

The first completeness proofs for interpretability logics \cite{jovelt90} went about this by labelling these roles inside the model. As such a single MCS could occur multiple times in a model with different labels. The labels that were employed in the old days were used to flag so-called \emph{criticality}. 

Criticality flagged that a specific MCS had a particular functionality in the Veltman model. In particular, for $C$ a formula, a $C$-critical successor flagged that the successor would avoid and not contain the formula $C$ in a strong sense so that this functionality propagates to parts accessible (either via $R$ or $S$) from that specific MCS. As such, completeness proofs could be very difficult and involved. Various different techniques were invented to keep track of all the different roles. 

Sometimes this could be done by keeping very close track of what roles could come after what other roles (\cite{jovelt90}, \cite{GorisJoosten:2012:SelfProvers}). In other occasions one had to consider many roles at the same time so as to avoid uncontrolled interaction between them (\cite{jonvelt99,Joosten:1998:MasterThesis,jogo08}). 

In 2004, the first author of this paper invented a slight variation of critical labels and called them \emph{assuring labels}. This variation now allowed to consider various roles at the same time. Where critical labels just keep track of a role with respect to one particular formula, the assuring labels actually are sets of formulas flagging that a MCS plays a role simultaneously with respect to all formulas in the set.

As a first application of assuring labels a one-page completeness proof of the logic \ilw was given in \cite{Bilkova-Goris-Joosten}. This should be contrasted with the very convoluted and complicated original completeness proof of around five pages based on criticality \cite{jonvelt99}. 

In years after the publication of \cite{Bilkova-Goris-Joosten}, almost all new  completeness proofs used the assuring labels
and in \cite{complgensem20} a uniform treatment of completeness proofs with respect to so-called generalised Veltman semantics could be given by making essential use of assuring labels.

The current paper is an expansion of \cite{Bilkova-Goris-Joosten} which was written on the occasion of Dick de Jongh's 60th birthday. A major draw-back of that paper is that it was actually written for Dick de Jongh and the paper assumed so much knowledge of the fields that virtually only Dick de Jongh could read it. Since the new technique has turned out to be so important, we decided to elaborate the old paper, make it self-contained, develop more of the theory and prove new results culminating in the current paper.

As such, Section \ref{sec:prelims} contains the needed technical preliminaries for the remainder of the paper. 
Then, in Section \ref{section:TechnicalIntroduction} we motivate the main notion of this paper: assuring labels. 

In Section \ref{sec:sassuring} we develop the general theory of assuring labels. Next, in Section \ref{sec:labellinglemmata} we shall see how assuring labels are good for imposing frame conditions on collections of MCSs. This was found to be useful in completeness proofs \cite{complgensem20}.

To illustrate the applicability, in Section \ref{sec:ilw} we give a short completeness and decidability proof of the logic \ilw and Section \ref{sec:finite} serves as a preparation.

Finally, Section \ref{sec:seq} analyses in a sense how iterations of labels may be needed in various situations, for example when considering the logic \extil{\principle{WR}}. As an illustration we prove completeness of \ilp with respect to a class of generalised Veltman frames where this iteration of labels is accounted for.

A shorter version of this preprint has been submitted as a journal paper (see
\cite{GorisBilkovaJoostenMikec:2022:JournalLabels}). The essential difference between the two versions is that this preprint also contains content from \cite{Bilkova-Goris-Joosten}, which has never been published as a paper and whose content fits naturally with our new findings.

\section{Preliminaries}\label{sec:prelims}
In this section we briefly revisit the main notions of interpretability logics that are relevant for this paper.

\subsection{Syntax and Logics for Interpretability}

The language of interpretability logics is given by 
\[
    \mathcal{F}::=\bot \mathbin{|} p \mathbin{|} \mathcal{F}\to \mathcal{F} \mathbin{|} \Box \mathcal{F} \mathbin{|} \mathcal{F}\rhd \mathcal{F} ,
\]
where $p$ ranges over a countable set of propositional variables. 
Other Boolean connectives are defined as abbreviations as usual. 
We treat $\rhd$ as having higher priority than $\to$, but lower than other logical connectives. 
We do not include $\Diamond$ in the language, rather we take $\Diamond A$ as an abbreviation for $\neg \Box \neg A$.

\begin{definition}
	The interpretability logic \il{} is axiomatised by the following axiom schemas.
	\begin{itemize}[leftmargin=2cm]
		\item[(Taut)] classical tautologies (in the new language);
		\item[(K)] $\square (A\rightarrow B)\rightarrow (\square A\rightarrow \square B)$;
		\item[(L)] $\square (\square A\rightarrow A)\rightarrow \square A$;
		\item[(J1)] $\square (A\rightarrow B)\rightarrow A\rhd B$;
		\item[(J2)] $(A\rhd B) \wedge (B\rhd C)  \rightarrow  A\rhd C$;
		\item[(J3)] $(A\rhd C) \wedge (B\rhd C)  \rightarrow  A\vee B\rhd C$;
		\item[(J4)] $A\rhd B\rightarrow(\Diamond A \rightarrow \Diamond B)$;
		\item[(J5)] $\Diamond A \rhd A$.
	\end{itemize}
	The rules of inference are Modus Ponens and Necessitation: $A/\Box A$.
\end{definition}

We will write $A\equiv B$ to denote $(A\rhd B) \wedge (B\rhd A)$. The following lemma is easy and we will use it throughout the paper, often tacitly.
Even though the proof is well-known and easy, we choose to include it as a warm-up for later reasoning so that we see the axioms at work.

\begin{lemma}\label{Theorem:BasicIlLemma}
The following are provable in \il.
\begin{enumerate}
\item \label{Theorem:BasicIlLemma:ItemTransitivity}
$\Box A\to \Box \Box  A$;

\item\label{Theorem:BasicIlLemma:ItemBoxVsRhd}
$\Box A \equiv \neg A \rhd \bot$;

\item\label{Theorem:BasicIlLemma:ItemLoebEquivClass}
$A\equiv A \wedge \Box \neg A$;

\item\label{Theorem:BasicIlLemma:ItemBoxFormulasVersusRhd}
$\Box C \wedge (A\rhd B) \ \to \ A \rhd B \wedge C$;

\item\label{Theorem:BasicIlLemma:ItemNoBoxDiamond}
For any formula $A$ we have
$\Box \Diamond A \equiv \Box \bot$.
\end{enumerate}
\end{lemma}

\begin{proof}
Item \ref{Theorem:BasicIlLemma:ItemTransitivity} is actually known to hold in \gl. We observe that in \il we can give an alternative proof: since $\Diamond \neg A \rhd \neg A$, by J4 we get $\Diamond \Diamond \neg A\to \Diamond \neg A$ and contraposition yields the required $\Box A\to \Box \Box  A$.

Item \ref{Theorem:BasicIlLemma:ItemBoxVsRhd} has two directions. First we observe that $\Box A \to \Box (\neg A \to \bot)$ so that by (J1) we obtain $\neg A \rhd \bot$. For the other direction, we apply (J4) to $\neg A \rhd \bot$ to obtain $\Diamond \neg A \to \Diamond \bot$. Since $\Diamond \bot$ is provably (actually in \gl) equivalent to $\bot$ we obtain $\Diamond \neg A \to \bot$ which is just $\neg \Diamond \neg A$, that is, $\Box A$.

Item \ref{Theorem:BasicIlLemma:ItemLoebEquivClass}. We will prove the first direction, and the other direction is proved similarly. We observe that $A \rhd (A\wedge \Box \neg A) \vee (A \wedge \Diamond A)$ so that by (J3) and (J2) we are done once we show $(A \wedge \Diamond A) \rhd A \wedge \Box \neg A$. By contraposing an instance of (L) we obtain $\Diamond A \to \Diamond (A \wedge \Box \neg A)$ so by Necessitation ($\Box(\Diamond A \to \Diamond (A \wedge \Box \neg A))$) and (J1) we obtain $\Diamond A \rhd  \Diamond (A \wedge \Box \neg A)$. Now (J5) yields $\Diamond (A \wedge \Box \neg A) \rhd A \wedge \Box \neg A$ so that (J2) gives $\Diamond A \rhd A \wedge \Box \neg A$. The result follows since clearly $A\wedge \Diamond A \rhd \Diamond A$.

Item \ref{Theorem:BasicIlLemma:ItemBoxFormulasVersusRhd} is easy since $\Box C \to \Box (B\to B \wedge C)$ so that $B\rhd B\wedge C$. Finally, Item \ref{Theorem:BasicIlLemma:ItemNoBoxDiamond} follows easily from $(L)$ since $\Box \Diamond A \to \Box \Diamond \top$ which implies $\Box (\Box \bot \to \bot)$ so that $\Box \bot$.
\end{proof}

By an \textit{extension} of \il we mean a logic which is (whose set of theorems is) a superset of (the set of theorems of) \il, and which is additionally closed under the same rules of inference. In this paper we shall consider extensions of \il with the following principles.
\[
\begin{array}{lll}
\principle{W} & \eqbydef & A\rhd B \rightarrow  A \rhd B \wedge \Box \neg A\\
\principle{M} & \eqbydef & A\rhd B \rightarrow  A \wedge \Box C \rhd B \wedge \Box C\\
\principle{P} & \eqbydef & A \rhd B \rightarrow \Box (A \rhd B) \\
\principle{M_0} & \eqbydef & A\rhd B \rightarrow \Diamond A \wedge \Box C \rhd B \wedge \Box C\\
\principle{R} & \eqbydef & A\rhd B \rightarrow  \neg (A \rhd \neg C) \rhd B\wedge \Box C\\
\end{array}
\]

Let us briefly comment on the significance of these principles. Montagna's principle \textsf{M}: $A \rhd B \to A \wedge \Box C \rhd B \wedge \Box C$ is valid in theories 
proving full induction. 
Berarducci \cite{bera90}  and Shavrukov \cite{shav88} independently proved that $\extil{}(T)=\extil{M}$ (\extil{} extended with \textsf{M}), if 
$T$ is $\Sigma_1$-sound and proves full induction, e.g.~\textbf{PA}.
The persistence principle \textsf{P}: $A \rhd B \to \Box(A \rhd B)$ is valid in finitely axiomatisable theories. 
Visser  \cite{Visser:1990:InterpretabilityLogic} proved the arithmetical completeness of \extil{P} w.r.t.\ any finitely 
axiomatisable $\Sigma_1$-sound theory 
containing $\mathsf I\Delta_0 + \mathsf{SUPEXP}$, where $\mathsf{SUPEXP}$ asserts the totality of 
superexponentiation, e.g.~\textbf{GB}. All other principles mentioned here are valid in all reasonable (sufficiently strong) theories (see \cite{Visser:1990:InterpretabilityLogic,jogo04,jogo08,jogo11,JoostenMasRoviraMikecVukovic:2021:OverviewGVSArXiv}).

\subsection{Relational semantics}

There are two basic and mutually related relational semantics\footnote{Apart from the various arithmetical semantics \cite{Visser:1991:FormalizationOfInterpretability,Joosten:2016:OreyHajek}, it is worth mentioning that recently topological semantics have been devised and studied too for interpretability logics and related logics \cite{iwata-kurahashi-2021}.} for interpretability logics. In both cases, the ternary relation $S$ will be conceived as a collection of parametrised binary relations. The first, and the most commonly used semantics, is \textit{Veltman semantics} (or \textit{ordinary Veltman semantics}). 
\begin{definition}
	A Veltman frame $\mathfrak{F}$ is a structure $(W,R, \{S_w : w\in W \})$, 
	where $W$ is a non-empty set, $R$ is a transitive and converse well-founded binary relation on 
	$W$ and for all $w\in W$ we have:
	\begin{itemize}
		\item[a)] $S_w\subseteq R[w]^2$, where $R[w]=\{x \in W : wRx\}$; 
		\item[b)] $S_w$ is  reflexive on $R[w]$; 
        \item[c)] $S_w$ is transitive;
		\item[d)] $S_w \supseteq R[w]^2 \cap R$, i.e.\ if $wRuRv$ then $uS_w v$.
	\end{itemize}
\end{definition} 

A \textit{Veltman model} is a quadruple $\mathfrak{M}=(W,R, \{S_w : w\in W \}, \Vdash)$,
    where the first three components form a Veltman frame.
The forcing relation $\Vdash$ is extended as usual in Boolean cases, 
and $w\Vdash A\rhd B$ holds if and only if for all $u$ such that $wRu$ and $u\Vdash A$ 
there exists $v$ such that $uS_w v$ and $v\Vdash B$.

The other commonly used semantics is the \textit{generalised Veltman semantics} by R.\ Verbrugge \cite{Verbrugge}.
\begin{definition}
	A generalised Veltman frame $\mathfrak F$ is a structure $(W,R,\{S_w:w\in W\})$, 
	where $W$ is a non-empty set, $R$ is a transitive and converse well-founded binary relation on 
$W$ and for all $w\in W$ we have:
	\begin{itemize}
		\item[a)] $S_w\subseteq R[w]\times \left(\mathcal{P}( {R[w]})\setminus\{\emptyset\}\right)$;
		\item[b)] $S_w$ is quasi-reflexive: $wRu$ implies $uS_w\{u\}$;
		\item[c)] $S_w$ is quasi-transitive: if $uS_wV$ and $vS_wZ_v$ for all $v\in V$, then  $uS_w(\bigcup_{v\in V}Z_v)$;
		\item[d)] if $wRuRv$, then $uS_w\{v\}$;
		\item[e)] monotonicity: if $uS_wV$ and $V\subseteq Z\subseteq R[w]$, then $uS_wZ$.
	\end{itemize}
\end{definition}
A \textit{generalised Veltman model} is a quadruple $\mathfrak{M}=(W,R, \{S_w : w\in W \}, \Vdash)$,
    where the first three components form a generalised Veltman frame. 
Sometimes it is useful to treat models as ordered pairs $(\mathfrak{F}, \Vdash)$ of a frame and a forcing relation.
With this semantics, $w\Vdash A\rhd B$ holds if and only if for all $u$ such that $wRu$ and $u\Vdash A$ 
there exists $V$ such that $uS_w V$ and $V\Vdash B$. 
By $V \Vdash B$ we mean $v \Vdash B$ for all $v \in V$.

For both semantics, we write $\mathfrak{M} \Vdash A$ (validity on a model) if $w \Vdash A$ for all $w \in W$, and $\mathfrak{F} \Vdash A$ (validity on a frame) if $(\mathfrak{F}, \Vdash) \Vdash A$ for all appropriate forcing relations $\Vdash$.

A \emph{frame condition} for an axiom scheme $A$ is a formula
$(A)$ (first or higher-order) in the language $\{ R, \{S_w : w\in W\} \}$ so that $\mathfrak{F}
\Vdash (A)$ (as a relational structure) if and only if $\langle
\mathfrak{F},V\rangle \Vdash A$ for every valuation~$V$.
For generalised semantics, we usually denote the frame condition by \kgen{A}.

Uppercase Greek, like $\Gamma$ and $\Delta$, will denote maximal consistent sets (MCSs).
It will be clear from the context with respect to what logic the consistency will refer.
Uppercase Roman denotes modal interpretability formulas $A,B,C,\ldots$ or sets
of such formulas $S,T,U,\ldots$.
An exception to this rule is that we might
write formulas from a set $S$ as $S_i$, $S_j$ etc.
In particular if $S$ is a set of formulas, then $\bigvee S_i$ denotes a finite
disjunction over some formulas in $S$.
If we talk of \textit{logics} or \textit{\extil{X}-theories} we mean extensions of \il that are closed under modus ponens and necessitation.  
As usual we use $\boxdot A$ as an abbreviation for $A\wedge\Box A$.
If $S$ is a set of formulas then we write $\Box S$ for $\{\Box A\mid A\in S\}$, and in general given some connective $\circ$, we write ${\circ} S$ for $\{ {\circ}A \mid A \in S \}$ and $S \circ Q$ for $\{ A \circ B \mid A \in S, B \in Q \}$. If $S$ and $Q$ are sets, we write $S \subseteq_{\sf fin} Q$ for ``$S$ is finite and $S \subseteq Q$''.

\section{Extending criticality}\label{section:TechnicalIntroduction}
In this section we motivate and define the notion of assuring successors and compare them to its restricted variant known as critical successors.

\subsection{Extending criticality}
As mentioned before, completeness proofs typically suit the following scheme. We take a formula $A$ that is not provable. Hence $\neg A$ is included in some MCS $\Gamma$. Next define the binary relation $R$ on MCSs together with the ternary relation $S$ so that the resulting structure is a model of the logic under consideration. Finally, we prove a so-called \emph{Truth lemma} that states
\begin{equation}\label{equation:TruthLemma}
    \forall \Delta \forall B \ \ \Big( B\in \Delta \ \Leftrightarrow \ \Delta \Vdash B \Big).
\end{equation}
Now, since $\neg A \in \Gamma$ we get that $A$ is falsified somewhere in our model.

It is easy to see that a least requirement for \eqref{equation:TruthLemma} to hold with respect to formulas of the form $\Box C$ is that whenever $\Gamma R \Delta$ we have for any $\Box A \in \Gamma$ that $A, \Box A\in \Delta$. This consideration gives rise to defining the following relation between MCSs.

\begin{definition}
For MCSs $\Gamma$ and $\Delta$ we define
\[
\Gamma \prec \Delta \ \ :\Leftrightarrow \ \ \forall A \ \Big( \Box A \in \Gamma \Rightarrow A,\Box A \in \Delta \Big).
\]
\end{definition}

We will now investigate what \eqref{equation:TruthLemma} imposes on the $S$ relation. In particular,. let us consider the condition for a formula $\neg (A \rhd B)$ to be true in some world $x$ in some particular model. From the previous section we know that $x\Vdash \neg (A\rhd B)$ if and only if there is some world $y$ so that $xRy$, so that $y\Vdash A$ but for no $z$ for which $yS_xz$ will we have $z\Vdash B$. In particular, since $yS_xy$ we see that $y\Vdash \neg B$. Moreover, since $yRu$ implies $yS_xu$ we also see that $y\Vdash \Box \neg B$.

Thus, certain transitions $\Gamma R \Delta$ actually should come with a promise that for any $\Delta'$ with $\Delta S_\Gamma \Delta'$ we will have $\neg B, \Box \neg B \in \Delta'$. Of course, we should also have $\neg C, \Box \neg C \in \Delta'$ for any $C$ so that $C\rhd B \in \Gamma$. 
Let us introduce the notion of criticality from \cite{jovelt90}.

\begin{definition}
For MCSs $\Gamma$ and $\Delta$ and for $C$ a formula, we say that $\Delta$ is a \emph{$C$-critical} successor of $\Gamma$ whenever
\[
\forall B \ \Big( B\rhd C \in \Gamma \ \Rightarrow \ \neg B, \Box \neg B \in \Delta \Big).
\]
We will write $\Gamma \prec^{C} \Delta$ in this case.\footnote{The usual notation for criticality is $\Gamma \prec_{C} \Delta$. We write $\Gamma \prec^{C} \Delta$ for criticality in this paper in order to more clearly distinguish it from assuringness, which we denote with $\Gamma \prec_{S} \Delta$.
We also note that the traditional notion of criticality, used e.g.~in \cite{jovelt90} and \cite{jonvelt99}, requires that there is some $\Box A \in \mathcal{D} \cap \Delta \setminus \Gamma$, where $\mathcal{D}$ is a finite set (usually a set of subformulas of some fixed formula). This additional requirement is a way of ensuring converse well-foundedness when criticality is used in the context of building Veltman models. We will not need this feature of criticality as we will use a different way of ensuring converse well-foundedness (also used in \cite{complgensem20}).
}
\end{definition}

It is easy to see that $C$-criticality naturally extends the $\prec$ relation as reflected by the following easy lemma.

\begin{lemma}
For MCSs $\Gamma$ and $\Delta$ we have $\Gamma \prec \Delta \ \Leftrightarrow \ \Gamma \prec^{\bot}\Delta$.
\end{lemma}
\begin{proof}
Immediate since $\Box A \in \Gamma \ \Leftrightarrow \ \neg A \rhd \bot \in \Gamma$.
\end{proof}
We can see $C$-criticality as a promise that the formula $C$ will be avoided in a strong sense. All completeness proofs before \cite{Bilkova-Goris-Joosten} made essential use of critical successors. Whenever in a structure of MCSs a $\Gamma\prec^C \Delta$ was there, the definition of the $S_\Gamma$ relation should reflect the promise that $C$ should be avoided. This strategy, although successful, resulted in a need for complicated book-keeping to keep all promises. 

An improvement can be made if we can deal with various promises at the same time. Suppose we wished to define $\Gamma \prec^{B,C} \Delta$ in such a way that it promises that both $B$ and $C$ are avoided in $\Delta$ in a strong sense. Requiring that simultaneously both $\Gamma \prec^{B} \Delta$ and $\Gamma \prec^{C} \Delta$ is not sufficient since the promises may interact. In particular 
\[
\mbox{if $A \rhd B\vee C \in \Gamma$ we should also require that $\neg A, \Box \neg A \in \Delta$.}
\]
It is this simple idea that adds a lot of power to the notion of criticality. However there is one more subtlety to it. It turns out to be very fruitful to apply a change of perspective. Instead of speaking of a promise to avoid certain formulas it turns out to be a very fruitful perspective to rather speak of \emph{assuring} certain formulas. If we do so, it will turn out that the set of promises has certain nice properties. In particular, it can be closed under logic consequence as proven in Lemma \ref{lemm:labiLukitsa}. These considerations give rise to the following definition.



\begin{definition}[Assuring successor]
\label{def:assuring_succ}
Let $S$ be a set of formulas, and $\Gamma$ and $\Delta$  some \extil{X}-MCSs (\extil{X} will usually be clear from context).
We define $\Gamma \sassuring{S} \Delta$, 
and say that $\Delta$ is an \emph{$S$-assuring successor} of $\Gamma$,
if for any formula $A$ and any finite $S'\subseteq S$ we have
$\neg A \rhd \bigvee_{S_j \in S'} \neg S_j \in \Gamma \Rightarrow A, \Box A \in \Delta$.
We will call $S$ a \emph{label} for $\Gamma$ and $\Delta$ or simply a \emph{label}.
\end{definition}

\begin{remark} \label{remark:even-more-positive-assuring}
The preceding definition is slightly different from all other definitions of an \emph{$S$-assuring successor} that have appeared so far (e.g.~\cite{Bilkova-Goris-Joosten} and \cite{complgensem20}). Compared to \cite{Bilkova-Goris-Joosten}, we do not require the existence of some $\Box C\in\Delta$ such that $\Box C\not\in\Gamma$ (this condition ensures converse well-foundedness in modal completeness proofs, but there are other way to ensure it, see e.g.~the definition of an \extil{X}-structure in \cite{complgensem20}). Compared to both papers, we require $\neg A \rhd \bigvee_{S_j \in S'} \neg S_j \in \Gamma \Rightarrow A, \Box A \in \Delta$ instead of $A \rhd \bigvee_{S_j \in S'} \neg S_j \in \Gamma \Rightarrow \neg  A, \Box \neg  A \in \Delta$. Here, the benefit is that the set of formulas which we will denote as $\Gamma_S^\boxdot$ will be closed under modus ponens and necessitation. This is the minimal set that satisfies the following property for all MCSs $\Delta$: $\Gamma_S^\boxdot \subseteq \Delta \Leftrightarrow \Gamma \prec_S \Delta$.
\end{remark}

In the following lemma we shall see that the notion of assuring successor on sets of formulas naturally extends the regular successor relation as well as the critical successor relation.

\begin{lemma}\label{lemma:assuringExtendsSuccessor}
\
\begin{enumerate}
\item
  $\Gamma \sassuring{\emptyset} \Delta \Leftrightarrow \Gamma \sucs \Delta$;

\item
  $\Delta$ is a $B$-critical successor of
  $\Gamma$ $\Leftrightarrow \Gamma \sassuring{\{ \neg B\}} \Delta$. 
\end{enumerate}
\end{lemma}

\begin{proof}
For the first item, we observe that the empty disjunction is per definition equivalent to $\bot$. Thus, by Lemma \ref{Theorem:BasicIlLemma}.\ref{Theorem:BasicIlLemma:ItemBoxVsRhd} we have $A\rhd \bot \in \Gamma$ if and only if $\Box \neg A \in \Gamma$.
Consequently, 
\[
\forall A\Big( \neg A\rhd \bot \in \Gamma \ \Rightarrow \  A, \Box  A \in \Delta  \Big) \ \Leftrightarrow \ \forall A\Big( \Box A \in \Gamma \ \Rightarrow \ A, \Box A \in \Delta  \Big).
\]

The $\Leftarrow$ direction of the second item is easy and the other direction follows from the first item of this lemma: if we take a finite subset of $\{ \neg B \}$ this is either the empty set, or $\{ \neg B \}$ itself. Now, $\neg A\rhd \neg \neg B \in \Gamma \Rightarrow A, \Box A \in \Delta$ follows from the assumption that $\Delta$ is a $B$-critical successor of $\Gamma$ and that $\Gamma$ is a MCS, and $\neg A\rhd \bot \in \Gamma \Rightarrow A, \Box A \in \Delta$ follows from the first item since critical successors are in particular successors.
\end{proof}

\subsection{Assuring labels assure}
As the name suggests, assuring labels assure certain formulas to be present. The relation \sassuring{S} assures elements in $\Gamma$ and $\Delta$ although it is not allowed to ``speak'' of consistency formulas, i.e.\ $\Diamond$-formulas cannot be contained in a label.
This is made explicit in the following lemma.

\begin{lemma} \label{lemma:assuringnessAssures}
The following holds:
\begin{enumerate}
\item
  $\Gamma \sassuring{S} \Delta \Rightarrow S , \Box S \subseteq \Delta$;
 
\item 
  $\Gamma \sassuring{S} \Delta \Rightarrow  
  \Diamond S \subseteq \Gamma$;

\item 
  $\Gamma \sassuring{S} \Delta \Rightarrow  
  \mbox{the label $S$ does not contain any formula of the form $\Diamond A$ or $\neg(A \rhd B)$.}$

\end{enumerate}
\end{lemma}

\begin{proof}
The first item is clear since for any $\sigma \in S$ we have that $\neg \sigma \rhd \neg \sigma$ is a theorem and whence in $\Gamma$. By the definition of $\Gamma \sassuring{S} \Delta$ we get that $\sigma, \Box \sigma \in \Delta$.

The second item follows from the first: since $\Gamma$ is maximal, for any $\sigma \in \Sigma$, either $\Diamond \sigma \in \Gamma$ or $\Box \neg \sigma \in \Gamma$. However, the latter would imply $\neg \sigma \in \Delta$ contradicting our first item.

For the last item we reason as follows. Suppose for a contradiction that there is some $\Diamond A$ in $S$. Then, by the first item we have both $\Diamond A \in \Delta$ and $\Box \Diamond A \in \Delta$. However, over \gl we have that $\Box \Diamond A$ is equivalent to $\Box \bot$. But $\Box \bot \in \Delta$ clearly contradicts $\Diamond A \in \Delta$. So, $\Diamond A \notin S$. This also implies $\neg(A \rhd B) \notin S$, since $\textbf{IL} \vdash \neg(A \rhd B) \to \Diamond A$ and $\Gamma$ is maximal.
\end{proof}

\section{The theory of assuring labels and of full labels}\label{sec:sassuring}
In this section we will expose theory of assuring labels and three special kind of assuring labels: full labels, maximal labels and complete labels.

\subsection{A general theory of assuring labels}
In the next section we will show how assuring successors can be used to solve, in a uniform way, certain problematic aspects of modal completeness proofs.

A label $S$ between $\Gamma \sassuring{S} \Delta$ keeps track of the formulas that are promised to be in 
$\Delta$ in virtue of certain interpretability formulas in $\Gamma$. The larger the label, the more promises it stores.

Often we can enlarge the label for free. To see how much we can add we need the following definition.

\newcommand{\boxdotset}[2]{{\ensuremath{{#2}^{\boxdot}_{#1}}}\xspace}
\begin{definition}\label{def:boxset}
For any set of formulas $T$ and maximal consistent set $\Delta$ we define
\begin{align*}
\boxset{T}{\Delta} &=
\{ \Box A \mid \textrm{ for some finite } T'\subseteq T, \neg A \rhd \bigvee_{T_i\in T'} \neg T_i \in \Delta \},\\
\boxdotset{T}{\Delta} &=
\{\Box A, A\mid \textrm{ for some finite } T'\subseteq T, \neg A\rhd\bigvee_{T_i\in T'}\neg T_i\in\Delta\}.
\end{align*}
\end{definition}

Our definition of these sets is slightly different than the standard definitions; see Remark \ref{remark:even-more-positive-assuring}. Note that $\boxset{\emptyset}{\Delta} =\{ \Box A \mid \neg A\rhd \bot \in \Delta \} = \{ \Box C \mid \Box C \in \Delta \}$. The latter equality holds because $\Delta$ is a MCS.
Furthermore, note that $\Gamma \prec_S \Delta$ holds precisely if $\boxdotset{S}{\Gamma} \subseteq \Delta$.
The next lemma tells us how promises propagate over composition of successors.

\begin{lemma} \label{lemm:trivialsassure}
For the relation \sassuring{S} we have the following observations.
\begin{enumerate}
\item
  $S\subseteq T \ \& \ \Gamma \sassuring{T} \Delta \Rightarrow \Gamma \sassuring{S} \Delta $;

\item\label{lemm:trivialsassure:itemSisPreserved}
  $\Gamma \sassuring{S} \Delta \sucs \Delta' \Rightarrow \Gamma \sassuring{S} \Delta'$.
\end{enumerate}
\end{lemma}

\begin{proof}
The first item is obvious since any finite subset of $S$ is also a finite subset of $T$ whenever $S\subseteq T$. For the second item we observe that $\Gamma\sassuring{S}\Delta$ implies $\boxdotset{S}{\Gamma} \subseteq \Delta$ whence by $\Delta \prec \Delta'$ and $\Box \boxdotset{S}{\Gamma} \subseteq \boxdotset{S}{\Gamma}$ we see that $\boxdotset{S}{\Gamma} \subseteq \Delta'$.
\end{proof}


\begin{notation}
Often we shall simply write $\bigvee \neg S_i$ to indicate some particular finite disjunction without really specifying it. If in the same context we will need another particular but otherwise unspecified disjunction we will flag this by using a different index. Thus $\bigvee \neg S_i \vee \bigvee \neg S_j$ stands for the disjunction of two particular but unspecified finite disjunctions of negated formulas from some label set $S$.

Often we will consider a finite collection of formulas $C_j$ such that each $C_j$ will interpret some finite disjunction of negated formulas from the label $S$. For each particular formula $C_j$ we will denote the corresponding disjunction by $\bigvee_k \neg S^j_k$ and thus write $C_j \rhd \bigvee_k \neg S^j_k$. Subsequently, we will denote the big disjunction over all $k$ and all corresponding $\neg S^k_j$ by $\bigvee_{k,j} \neg S^k_j$ so that $\bigvee C_k \rhd \bigvee_{k,j} \neg S^k_j$. 

Often we shall further relax notation by omitting these sub-indices so that, par abus de langage, we will write $\bigvee \neg S^k_j$ where the context will make clear if one should read $\bigvee_{k} \neg S^k_j$ or $\bigvee_{k,j} \neg S^k_j$.
\end{notation}

The following lemma gives us a way to extend labels.

\begin{lemma}\label{lemm:labil}
  For any logic \textup{(}i.e.\ extension of $\extil{}$\textup{)} we have\footnote{
  In this type of claims, the choice of logic dictates what type of maximality we can assume for $\Gamma$, $\Delta$ and any other set appearing in the statement. If we are not explicit about it, we require only that $\Gamma$, $\Delta$ etc. are \extil{}-MCSs.
  }
  $\Gamma\sassuring{S}\Delta\Rightarrow\Gamma\sassuring{S\cup\boxdotset{S}{\Gamma}}\Delta$.
\end{lemma}

\begin{proof}
  Suppose $\Gamma\sassuring{S}\Delta$ and
  $\neg C\rhd\bigvee_i \neg S_i\vee\bigvee_j (\neg A_j\vee\Diamond \neg A_j)\in\Gamma$ for some finite collections of formulas $S_i \in S$ and $A_j, \Box A_j \in \boxdotset{S}{\Gamma}$. In particular, for each $j$ we have $\neg A_j \rhd \bigvee_k \neg S_k^j \in \Gamma$ for some finite collection (depending on $j$) of formulas $S_k^j$ from the label $S$.
  Then, since for each $A_j$ we have $\Diamond \neg A_j \rhd \neg A_j$ we obtain $\neg C\rhd\bigvee_i\neg S_i\vee\bigvee_j \neg A_j\in\Gamma$ and thus
  $\neg C\rhd\bigvee_i\neg S_i\vee\bigvee_{j,k} \neg S^j_k\in\Gamma$ which implies
  $C,\Box C\in\Delta$ since we assumed $\Gamma\sassuring{S}\Delta$. 
\end{proof}

This lemma tells us in a sense that when we have $\Gamma \sassuring{S}\Delta$, then certain sentences in $\Gamma$ justify that we may extend the label $S$. Will likewise the occurrence of sentences in $\Delta$ allow us to extend the label $S$? The next lemma tells us that this is not the case. In particular, if $\neg A\rhd \bigvee \neg S_i$ for some $S_i \in S' \subseteq_{\sf fin}S$, then by definition $A, \Box A \in \Delta$. However, when for some arbitrary $A$ we have $A, \Box A \in \Delta$, this does not allow us to extend our label $S$.

\begin{lemma}
\label{lemma:box_vs_assuringness}
There is a $\Gamma \sassuring{S}\Delta$ \textup{(}$\Gamma$ and $\Delta$ are \extil{}-MCSs\textup{)} with $p, \Box p \in \Delta$ but $\neg(\Gamma \sassuring{S\cup \{ p\}}\Delta)$.
\end{lemma}

\begin{proof}
\begin{figure}
    \centering
    \scalebox{0.8}{
        \tikzset{every picture/.style={line width=0.75pt}} 
        
        \begin{tikzpicture}[x=0.75pt,y=0.75pt,yscale=-1,xscale=1]
        
        \draw    (380.5,543) -- (380.5,484) ;
        \draw [shift={(380.5,482)}, rotate = 450] [color={rgb, 255:red, 0; green, 0; blue, 0 }  ][line width=0.75]    (10.93,-3.29) .. controls (6.95,-1.4) and (3.31,-0.3) .. (0,0) .. controls (3.31,0.3) and (6.95,1.4) .. (10.93,3.29)   ;
        \draw    (380.5,482) .. controls (420.1,452.3) and (434.22,511.79) .. (473.31,483.87) ;
        \draw [shift={(474.5,483)}, rotate = 503.13] [color={rgb, 255:red, 0; green, 0; blue, 0 }  ][line width=0.75]    (10.93,-3.29) .. controls (6.95,-1.4) and (3.31,-0.3) .. (0,0) .. controls (3.31,0.3) and (6.95,1.4) .. (10.93,3.29)   ;
        \draw    (380.5,543) -- (472.81,484.08) ;
        \draw [shift={(474.5,483)}, rotate = 507.45] [color={rgb, 255:red, 0; green, 0; blue, 0 }  ][line width=0.75]    (10.93,-3.29) .. controls (6.95,-1.4) and (3.31,-0.3) .. (0,0) .. controls (3.31,0.3) and (6.95,1.4) .. (10.93,3.29)   ;
        
        \draw (381,549.67) node   [align=left] {$\displaystyle x$};
        \draw (373,475.67) node   [align=left] {$\displaystyle y$};
        \draw (483,476.67) node   [align=left] {$\displaystyle z$};
        \draw (370.36,515) node   [align=left] {$\displaystyle S$};
        \draw (378,459) node   [align=left] {$\displaystyle q,\ p,\ \square p$};
        \draw (484,460) node   [align=left] {$\displaystyle \neg p$};
        \draw (388,562) node   [align=left] {$\displaystyle q\rhd \neg p$};

        \end{tikzpicture}
    }
    \caption{Situation described in Lemma \ref{lemma:box_vs_assuringness}.}
    \label{fig:box_vs_assuringness}
\end{figure}
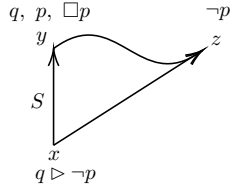

Consider the model consisting of three points $x,y$ and $z$ given in Figure \ref{fig:box_vs_assuringness}. We take $\Gamma$ to be the modal theory of $x$ and $\Delta$ to be the modal theory of $y$. Since $q\in \Delta$ and $(q\rhd \neg p)\in \Gamma$, whatever we take for $S$ with $\Gamma \sassuring{S}\Delta$, we will never have $\Gamma \sassuring{S\cup \{ p\}}\Delta$.
\end{proof}

\subsection{Full labels}
In Lemma \ref{lemm:labil} we saw that the $S$-assuringness between two sets $\Gamma \sassuring{S}\Delta$ can only be automatically extended via $\Gamma$. The next lemma tells us that there are other ways to `freely extend' a label.

\begin{lemma}\label{lemm:labiLukitsa}
  For any logic we have
  \begin{enumerate}
      \item 
  $(\Gamma\sassuring{S}\Delta) \ \mathbin{\&} \  (S\vdash \varphi)\ \  \Longrightarrow\ \ \Gamma\sassuring{S\cup \{  \varphi\}}\Delta$;
  
      \item 
  $\Gamma\sassuring{S}\Delta \ \  \Longrightarrow\ \ \Gamma\sassuring{S\cup \Box S}\Delta$;
  
  \end{enumerate}
\end{lemma}

\begin{proof}
For the first item---that labels can be closed under logical consequence---we assume that $S\vdash \varphi$ where the notion of logical consequence depends on the logic in question. Thus for some $S_1, \ldots, S_n \in S$ we have $S_1 \wedge \ldots \wedge S_n \vdash \varphi$. Consequently, $\vdash \neg \varphi \to \bigvee \neg S_j$ and also $\vdash \Box (\neg \varphi \to \bigvee \neg S_j)$. Thus, if $\Gamma \sassuring{S}\Delta$ and $(\neg A \rhd \bigvee \neg S_i \vee \neg \varphi) \in \Gamma$, also $(\neg A \rhd \bigvee \neg S_i) \in \Gamma$ so that $A, \Box A \in \Delta$ and we conclude $\Gamma \sassuring{S \cup \{  \varphi \}}\Delta$. 

For the second item, we consider $(\neg A\rhd \bigvee \neg S_i \vee \bigvee \neg \Box S_j) \in \Gamma$. But since $\neg \Box S_j \equiv \Diamond \neg S_j$ and $\Diamond \neg S_j \rhd \neg S_j$ we conclude $(\neg A\rhd \bigvee \neg S_i \vee \bigvee \neg S_j)$ so that $A, \Box A \in \Delta$.
\end{proof}

This Lemma \ref{lemm:labiLukitsa} tells us that given an extension $\extil{X}$ of the logic \extil{}, we can freely extend labels to be closed under $\extil{X}$-consequences and to be closed under necessitation. Thus, we can identify labels with \ilx-theories.  

Moreover, Lemma \ref{lemm:labil} tells us that we can freely close off a label $S$ for  $\Gamma\sassuring{S}\Delta$ under $\boxdotset{S}{\Gamma}$.
These observations lead us to the definition of \emph{$\Gamma$-full labels}. When the context makes clear which $\Gamma$ is meant we shall simply speak of \emph{full labels}.

\begin{definition}\label{definition:fullLabel}
For $\Gamma$ an \extil{X}-MCS we call a set $S$ of formulas  a \emph{$\Gamma$-full label} whenever we have the following:
\begin{enumerate}
    
    \item \label{definition:fullLabel:BoxDotIn}
    $\neg A\rhd \bigvee \neg S_i \in \Gamma \ \Longrightarrow \ A , \Box A \in S$; 
    
    
    
    
    \item\label{definition:fullLabel:LogicalConsequence}
    the label $S$ is closed under logical consequence (w.r.t.~\extil{X}), that is, if $S\vdash \varphi$, then $\varphi \in S$;
    
    \item\label{definition:fullLabel:Necessitation}
    the label $S$ is closed under necessitation, that is, if $B \in S$, then $\Box B \in S$.
    
\end{enumerate}
\end{definition}

A direct consequence of Item \ref{definition:fullLabel:BoxDotIn} of this definition is that in particular $\Box A \in \Gamma$ implies $A\in S$.
The following lemma is easy and could be taken as an alternative definition.

\begin{lemma}\label{theorem:alternativeDefinitionFullLabel}
Given a MCS $\Gamma$ and a set of formulas $S$, then $S$ is a \emph{$\Gamma$-full label} if and only if $S$ is an $\extil{X}$-theory extending $\boxdotset{S}{\Gamma}$.
\end{lemma}

\begin{proof}
Both directions are easy. We observe that Items \ref{definition:fullLabel:LogicalConsequence} and \ref{definition:fullLabel:Necessitation} of Definition \ref{definition:fullLabel} exactly state that $S$ is an \ilx-theory whereas Item \ref{definition:fullLabel:BoxDotIn} corresponds to $\boxdotset{S}{\Gamma} \subseteq S$.  
\end{proof}

If we stick to full labels, there is a close correspondence between theories and labels. We find this observation so essential that we formulate it explicitly as a lemma:

\begin{lemma}
For any logic \extil{X}, if $\Gamma\sassuring{S}\Delta$ and $S$ is full, then $S$ is a consistent \extil{X}-theory.
\end{lemma}

A natural question readily suggests itself: which \ilx-theories may occur as label between two \ilx-MCSs? Lemma \ref{theorem:sufficientConditionForATheoryToBeLabel} below tells us that it is a sufficient condition that $S \nvdash_{\ilx} \Box \bot$.

Full labels contain as many free promises as possible and posses certain nice closure properties. In particular, we have the following lemma that justify the name `full'.

\begin{lemma}
\label{lemma:nice_characterisation}
Given an \extil{X}-\mcs $\Gamma$ and a label $S$, the label $S$ is \emph{$\Gamma$-full} if and only if 
\[
\forall T \Big( S\subseteq T \wedge \forall \Delta \big( \Gamma\sassuring{S}\Delta \Rightarrow \Gamma\sassuring{T}\Delta\big) \ \Longrightarrow \ S = T  \Big).
\]
The $S$ and $T$ range here over $\extil{X}$-theories, and $\Delta$ over \extil{X}-MCSs. 
\end{lemma}

\begin{proof}
First assume that $S$ is a $\Gamma$-full label and $S\subset T$. We want to show there is a \mcs $\Delta$ with $\Gamma \sassuring{S}\Delta$ but $\neg (\Gamma \sassuring{T}\Delta)$. As $S\subset T$, there is some $\phi\in T$ for which we have $\phi\notin S$, and therefore, by $S$ being $\Gamma$-full, $\phi\notin \boxdotset{S}{\Gamma}\subseteq S$. Since $S$ is a theory also $S\nvdash \phi$ and $\boxdotset{S}{\Gamma}\nvdash \phi$. Then there exists a \mcs $\Delta$ containing $\boxdotset{S}{\Gamma}$ with $\phi\notin\Delta$. Clearly $\Gamma \sassuring{S}\Delta$, and as $\neg\phi\rhd\neg\phi\in\Gamma, \phi\in T$ and $\phi\notin\Delta$, we see that $\neg (\Gamma \sassuring{T}\Delta)$. 

For the other direction assume $\boxdotset{S}{\Gamma}\not\subseteq S$. We want to find a theory $T\supset S$ with $\forall \Delta \big( \Gamma\sassuring{S}\Delta \Rightarrow \Gamma\sassuring{T}\Delta\big)$. Take $T$ to be the theory generated by $S\cup \boxdotset{S}{\Gamma}$: since $\boxdotset{S}{\Gamma}\not\subseteq S$, it certainly is bigger than $S$. Now assume $\Gamma\sassuring{S}\Delta$, but then $\Gamma\sassuring{S\cup \boxdotset{S}{\Gamma}}\Delta$ by Lemma \ref{lemm:labil}.
\end{proof}

We shall now see that each label can be extended to a full label. To this end, let us first observe that (full) labels are closed under intersections.

\begin{lemma}\label{theorem:labelsClosedUnderIntersections}
Let $\Gamma$ and $\Delta$ be \ilx-MCSs and let $\{S_i \}_{i\in I}$ be a collection of labels so that $\Gamma \sassuring{S_i} \Delta$. We have that
\[
\Gamma \sassuring{\bigcap_{i\in I}S_i} \Delta.
\]
Furthermore, if each $S_i$ is full, then so is $\bigcap_{i\in I}S_i$.
\end{lemma}

\begin{proof}
That $\Gamma \sassuring{\bigcap_{i\in I}S_i} \Delta$ is immediate by the first item of Lemma \ref{lemm:trivialsassure} and it is easy to check that fullness is preserved under taking intersections.
\end{proof}

We can see the closure properties of full sets as defining monotone operators so that clearly, each label is included in some full label. Given this and the previous lemma it now makes sense to define the closure of a label as the intersection of all the full labels containing it. 

\begin{definition}
For given \ilx-MCSs $\Gamma$ and $\Delta$ with $\Gamma \sassuring{S} \Delta$ we define the closure $\overline{S(\Gamma)}$ of $S$ to be the smallest $\Gamma$-full label that contains $S$ so that $\Gamma \sassuring{\overline{S(\Gamma)}} \Delta$.
\end{definition}


The precise nature of $\overline{S(\Gamma)}$ depends on $\Gamma$ but when $\Gamma$ is clear from the context, we shall simply write $\overline{S}$.
The next lemma expresses that the closure of a label has all the properties that one would expect from a closure operator.

\begin{lemma}\label{theorem:fullClosure}
Let $\Gamma, \Delta$ be \ilx-MCSs and $S$ be a set of formulas. We have
\begin{enumerate}

    \item If $\Gamma \sassuring{S} \Delta$, then $\overline{S} = \Gamma^\boxdot_S$;
    
    \item
$\Gamma \sassuring{S} \Delta \ \ \Leftrightarrow \ \ \Gamma \sassuring{\overline{S}} \Delta$;
\\\\
For the remaining items we assume $\Gamma \sassuring{S} \Delta$.
\item 
    $\overline{S}$ is $\Gamma$-full;

\item
$\overline{S} = \overline{\overline{S}}$;

\item
$\overline{S}$ is the smallest \ilx-theory containing $S\cup \Gamma^\boxdot_S$.

\end{enumerate}

\end{lemma}

\begin{proof}

To prove that $\overline{S} = \Gamma^\boxdot_S$, we need to show that $\Gamma^\boxdot_S$ is contained in all $\Gamma$-full labels that extends $S$ (which is a simple consequence of Lemma \ref{theorem:alternativeDefinitionFullLabel}), and that it is a $\Gamma$-full label itself (which remains to be shown).
The latter is (again by Lemma \ref{theorem:alternativeDefinitionFullLabel}) equivalent to the claim that $\Gamma^\boxdot_S$ is a logic and  $\Gamma_{ \Gamma^\boxdot_S  }^\boxdot \subseteq \Gamma^\boxdot_S$. 

Let us first show that $\Gamma^\boxdot_S$ is a logic. 
Let $A, A \to B \in \Gamma^\boxdot_S$. It follows (possibly using \textbf{J5}) that $\neg A \rhd \bigvee_i \neg S_i \in \Gamma$ and $\neg(A \to B) \rhd \bigvee_j \neg S_j \in \Gamma$, so $\neg B \rhd \neg A \vee \neg(A \to B) \rhd \bigvee_i \neg S_i \vee \bigvee_j \neg S_j \in \Gamma$. Thus, $B \in \Gamma^\boxdot_S$. Next we verify that $\Gamma^\boxdot_S$ is closed under necessitation. Assume $A, \Box A \in \Gamma^\boxdot_S$ due to $\neg A \rhd \bigvee \neg S_i \in \Gamma$ for some $S_i \in S$. We are to show $\Box A, \Box \Box A \in \Gamma^\boxdot_S$. Notice that $\Box A \in \Gamma^\boxdot_S$ is a part of our assumption, so we only need to show $\Box \Box A \in \Gamma^\boxdot_S$. Clearly $\neg \Box \Box A \rhd \neg A \rhd \bigvee \neg S_i \in \Gamma$, so $\Box \Box A \in \Gamma^\boxdot_S$.

Let us show that $\Gamma_{ \Gamma^\boxdot_S }^\boxdot \subseteq \Gamma^\boxdot_S$.
Assume $A, \Box A \in \Gamma_{ \Gamma^\boxdot_S }^\boxdot$ due to $\neg A \rhd \bigvee \neg B_i \in \Gamma$ for some $B_i \in \Gamma^\boxdot_S$. 
For each formula $B_i$ there are some formulas $S_j^i \in S$ such that $\neg B_i \rhd \bigvee \neg S_j^i \in \Gamma$. Applying the axiom schema \textbf{J3} a certain number of times, we get $\neg A \rhd \bigvee_{i,j} \neg S_j^i \in \Gamma$, thus $A, \Box A \in \Gamma^\boxdot_S$, as required.

For the second claim, by Lemma \ref{lemm:labil} we have that $\Gamma \sassuring{S} \Delta$ implies $\Gamma \sassuring{ \boxdotset{S}{\Gamma} } \Delta$. The other direction is simple as $S \subseteq \boxdotset{S}{\Gamma}$. 
The third claim holds by definition. 
For the forth claim, notice that as $\overline{S}$ is $\Gamma$-full, it itself is its smallest $\Gamma$-full extension. 
For the fifth claim, notice that $\overline{S}$ is an \extil{X}-theory (because it is a $\Gamma$-full label, and also because by the first claim it equals $\boxdotset{S}{\Gamma}$ which is a theory), and furthermore $\overline{S} = \boxdotset{S}{\Gamma} = S \cup \boxdotset{S}{\Gamma}$.

\end{proof}



Full labels can at times simplify matters. In particular, they clearly propagate along successors as expressed by the following lemma.

\begin{lemma}\label{theorem:informationAccrualOverFullLabels}
For any logic, if $\Gamma \sassuring{S}\Delta \sassuring{T}\Lambda$ for some full labels $S$ and $T$, then $S\subseteq T$.
\end{lemma}

\begin{proof}
For any $S_i\in S$ we have $\Box S_i \in \Delta$ so, by fullness, $S_i \in T$.
\end{proof}

Thus, this Lemma states that full labels accrue information along the top successor relation. Does information between related full labels also `reflect down'? To put it otherwise, it may be natural to ask if Lemma \ref{lemm:trivialsassure}.\ref{lemm:trivialsassure:itemSisPreserved} ($\Gamma \sassuring{S} \Delta \sucs \Delta' \Rightarrow \Gamma \sassuring{S} \Delta'$) can be strengthened. That is to say, suppose we have $\Gamma \sassuring{S} \Delta \sassuring{T} \Delta'$, can we say something more than just $\Gamma \sassuring{S} \Delta'$? As we shall see in the next section, it turns out that for extensions of \il we often can. 
In general this does not seem to hold, at least if we do not require our labels to be full. 
Suppose $\Gamma \sassuring{\emptyset}\Delta\sassuring{ \{ p \} }\Delta'$ (see Figure \ref{fig:downward}).
If $p \rhd \neg p \in \Gamma$ and $p \in \Delta'$,
    there is a MCS $\Lambda$ with $\neg p \in \Lambda$.
Clearly, the fact that we have $\{ p \}$ between $\Delta$ and $\Delta'$ did not stop $\neg p \in \Lambda$.
Let us mention a question that we do not have a definite answer for. 
Suppose $\Gamma \sassuring{S}\Delta\sassuring{T}\Delta'$ and $S$ and $T$ are full labels.
Is there a (non-trivial) notion of a ``$T$-influenced formula'' such that we may put the $T$-influenced formulas between $\Gamma$ and $\Delta'$?

\begin{figure}
    \centering
    \begin{minipage}[b]{0.45\textwidth} 
        \centering
        \scalebox{0.8}{
            \begin{tikzpicture}[x=0.75pt,y=0.75pt,yscale=-1,xscale=1]
            
            \draw    (134.5,567) -- (134.5,508) ;
            \draw [shift={(134.5,506)}, rotate = 450] [color={rgb, 255:red, 0; green, 0; blue, 0 }  ][line width=0.75]    (10.93,-3.29) .. controls (6.95,-1.4) and (3.31,-0.3) .. (0,0) .. controls (3.31,0.3) and (6.95,1.4) .. (10.93,3.29)   ;
            \draw    (134.5,491) -- (134.5,432) ;
            \draw [shift={(134.5,430)}, rotate = 450] [color={rgb, 255:red, 0; green, 0; blue, 0 }  ][line width=0.75]    (10.93,-3.29) .. controls (6.95,-1.4) and (3.31,-0.3) .. (0,0) .. controls (3.31,0.3) and (6.95,1.4) .. (10.93,3.29)   ;
            \draw    (134.5,430) .. controls (174.1,400.3) and (188.22,459.79) .. (227.31,431.87) ;
            \draw [shift={(228.5,431)}, rotate = 503.13] [color={rgb, 255:red, 0; green, 0; blue, 0 }  ][line width=0.75]    (10.93,-3.29) .. controls (6.95,-1.4) and (3.31,-0.3) .. (0,0) .. controls (3.31,0.3) and (6.95,1.4) .. (10.93,3.29)   ;
            \draw    (134.5,567) -- (227.36,432.65) ;
            \draw [shift={(228.5,431)}, rotate = 484.65] [color={rgb, 255:red, 0; green, 0; blue, 0 }  ][line width=0.75]    (10.93,-3.29) .. controls (6.95,-1.4) and (3.31,-0.3) .. (0,0) .. controls (3.31,0.3) and (6.95,1.4) .. (10.93,3.29)   ;
            
            \draw (135,573.67) node   [align=left] {$\displaystyle \Gamma$};
            \draw (135,497.67) node   [align=left] {$\displaystyle \Delta$};
            \draw (127,424.67) node   [align=left] {$\displaystyle \Delta'$};
            \draw (237,424.67) node   [align=left] {$\displaystyle \ \ \ \Lambda \ni \neg p$};
            \draw (179,440.67) node   [align=left] {$\displaystyle \Gamma$};
            \draw (124.36,467) node   [align=left] {$\displaystyle \{ p \}$};
            \draw (123.36,537) node   [align=left] {$\displaystyle \emptyset $};
            \draw (191.36,506) node   [align=left] {$\displaystyle \emptyset $};
            \end{tikzpicture}
        }
        \caption{Downward influence}
        \label{fig:downward}
    \end{minipage}
    \begin{minipage}[b]{0.45\textwidth} 
        \centering
        \scalebox{0.8}{
        \begin{tikzpicture}[x=0.75pt,y=0.75pt,yscale=-1,xscale=1]
        
        \draw    (484.5,193.39) -- (423.17,107.92) ;
        \draw [shift={(422,106.3)}, rotate = 414.34000000000003] [color={rgb, 255:red, 0; green, 0; blue, 0 }  ][line width=0.75]    (10.93,-3.29) .. controls (6.95,-1.4) and (3.31,-0.3) .. (0,0) .. controls (3.31,0.3) and (6.95,1.4) .. (10.93,3.29)   ;
        \draw    (484.5,193.39) -- (360.14,106.48) ;
        \draw [shift={(358.5,105.33)}, rotate = 394.95] [color={rgb, 255:red, 0; green, 0; blue, 0 }  ][line width=0.75]    (10.93,-3.29) .. controls (6.95,-1.4) and (3.31,-0.3) .. (0,0) .. controls (3.31,0.3) and (6.95,1.4) .. (10.93,3.29)   ;
        \draw    (422,106) .. controls (383.09,91.55) and (396.09,122.38) .. (360.18,106.11) ;
        \draw [shift={(358.5,105.33)}, rotate = 385.35] [color={rgb, 255:red, 0; green, 0; blue, 0 }  ][line width=0.75]    (10.93,-3.29) .. controls (6.95,-1.4) and (3.31,-0.3) .. (0,0) .. controls (3.31,0.3) and (6.95,1.4) .. (10.93,3.29)   ;
        \draw    (484,193.39) -- (545.33,107.92) ;
        \draw [shift={(546.5,106.3)}, rotate = 485.66] [color={rgb, 255:red, 0; green, 0; blue, 0 }  ][line width=0.75]    (10.93,-3.29) .. controls (6.95,-1.4) and (3.31,-0.3) .. (0,0) .. controls (3.31,0.3) and (6.95,1.4) .. (10.93,3.29)   ;
        \draw    (484,193.39) -- (608.36,106.48) ;
        \draw [shift={(610,105.33)}, rotate = 505.05] [color={rgb, 255:red, 0; green, 0; blue, 0 }  ][line width=0.75]    (10.93,-3.29) .. controls (6.95,-1.4) and (3.31,-0.3) .. (0,0) .. controls (3.31,0.3) and (6.95,1.4) .. (10.93,3.29)   ;
        \draw    (546.5,106) .. controls (585.41,91.55) and (572.41,122.38) .. (608.32,106.11) ;
        \draw [shift={(610,105.33)}, rotate = 514.65] [color={rgb, 255:red, 0; green, 0; blue, 0 }  ][line width=0.75]    (10.93,-3.29) .. controls (6.95,-1.4) and (3.31,-0.3) .. (0,0) .. controls (3.31,0.3) and (6.95,1.4) .. (10.93,3.29)   ;
        
        \draw (436,86) node   [align=left] {$\displaystyle p,q,r$};
        \draw (438.5,99.33) node   [align=left] {$\displaystyle u_{1}$};
        \draw (352.5,100.33) node   [align=left] {$\displaystyle v_{1}$};
        \draw (354,86) node   [align=left] {$\displaystyle p$};
        \draw (533,86) node   [align=left] {$\displaystyle p,q,r$};
        \draw (535.5,99.33) node   [align=left] {$\displaystyle u_{2}$};
        \draw (621.5,100.33) node   [align=left] {$\displaystyle v_{2}$};
        \draw (623,86) node   [align=left] {$\displaystyle q$};
        \draw (503,213.45) node   [align=left] {$\displaystyle \ r\rhd \neg p\lor \neg q,\ \neg ( r\rhd \neg p) ,\ \neg ( r\rhd \neg q)$};
        \draw (461,136) node   [align=left] {$\displaystyle \{p\}$};
        \draw (507,136) node   [align=left] {$\displaystyle \{q\}$};
        \draw (484.5,200.33) node   [align=left] {$\displaystyle w$};
        \end{tikzpicture}
        }
        
        \caption{Incomparable labels}
        \label{fig:no_maximum_ex}
    \end{minipage}
\end{figure}

\subsection{Maximal and complete labels}

Although a label can be full, this does not mean we can always find a maximum among the possible labels. We shall now exhibit a model that generates maximal consistent sets $\Gamma$ and $\Delta$ with two incomparable labels between them. Thus, full labels need not necessarily be maximum.

\begin{lemma}\label{theorem:labelsNotClosedUnderUnions}
There are \extil{}-MCSs $\Gamma$ and $\Delta$ and labels $S$ and $T$ with $\Gamma \sassuring{S}\Delta$ and $\Gamma \sassuring{T}\Delta$ so that $\neg(\Gamma \sassuring{S\cup T}\Delta)$.
\end{lemma}

\begin{proof}
We let $S:= \{p\}$, $T:= \{ q \}$ and consider the model in Figure \ref{fig:no_maximum_ex}. Let $\Gamma$ be the modal theory of the world $w$ and $\Delta$ be the modal theory of the world $u_1$. Clearly, $u_1$ and $u_2$ have the same modal theory, and by modal soundness these sets are (maximal) consistent sets.
\end{proof}

The above lemma exhibits two incomparable labels and shows that in general we cannot expect there to be a label between two MCSs that is a maximum with respect to the inclusion relation. However, we can always find a maximal label as we shall observe now. 

In Lemma \ref{theorem:labelsClosedUnderIntersections} we learned that full labels are always closed under intersections. In Lemma \ref{theorem:labelsNotClosedUnderUnions} we saw that labels are not necessarily closed under unions. However, we now observe that (full) labels are also closed under unions of chains.

\begin{lemma}\label{theorem:labelsClosedUnderUnions}
Let $\Gamma$ and $\Delta$ be \ilx-MCSs and let $\{ S_i\}_{i\in I}$ be a chain of labels with $\Gamma \sassuring{S_i} \Delta$ and $S_i \subset S_{i+1}$. We have that
\[
\Gamma \sassuring{\bigcup_{i\in I}S_i} \Delta.
\]
Moreover, if each of the $S_i$ are $\Gamma$ full, then so is $\bigcup_{i\in I}S_i$.
\end{lemma}

\begin{proof}
The lemma follows directly from the definition of assuringness since this involves only finitely many formulas from the label.
\end{proof}

As usual, closure under unions of chains of labels gives us the existence of a maximal label.

\begin{lemma}
Let $\Gamma$ and $\Delta$ be \ilx-MCSs with $\Gamma \sassuring{S} \Delta$ for some \ilx-theory $S$. Then, we can find a maximal \ilx-theory $S'$ extending $S$ so that $\Gamma \sassuring{S'} \Delta$.
\end{lemma}

\begin{proof}
By the previous lemma and Zorn's Lemma.
\end{proof}

It is easy to see that each maximal label is itself a full label. The next easy lemma tells us that indeed there are full labels that are not maximal.

\begin{lemma}
There exist MCSs $\Gamma, \Delta$ and $\Gamma$-full labels $S, T$ with $S \subsetneq T$ so that $\Gamma \sassuring{T} \Delta$.
\end{lemma}

\begin{proof}
We consider the model with domain $\{ x, y, y' \}$ so that all the $R$ relations are given by $xRy$ and $xRy'$ and $S_x$ is minimal ($yS_x y$ and $y'S_x y'$). Furthermore, $y\Vdash p,q$ and no other propositional variables and  $y'\Vdash p$ and no other propositional variables. We denote the modal theory of a world by its corresponding uppercase letter. Clearly $X \sassuring{\{p,q\}} Y$ whence by Lemma \ref{theorem:fullClosure}, also $X \sassuring{\overline{\{p,q\}}} Y$. Of course we also have $X \sassuring{\{p\}} Y$ whence $X \sassuring{\overline{\{p\}}} Y$. Now $\overline{\{p\}}$ cannot be equal to $\overline{\{p,q\}}$ since $X \sassuring{\{ p\}} Y'$ whence 
$X \sassuring{\overline{\{p\}}}$ and since $q\notin Y'$ we cannot have 
$X \sassuring{\overline{\{p,q\}}} Y'$ whence $\overline{\{p\}}\subsetneq \overline{\{p,q\}}$. Thus we have proven the lemma with $\Gamma = X$, $\Delta = Y$, $S= \overline{\{p\}}$ and $T= \overline{\{p,q\}}$.
\end{proof}

It is clear that not all assuring labels are full. Thus, the three notions of assuring, full and maximal labels define labels of increasing information. We can add an additional increment in information by considering the following definition.

\begin{definition}
We call a $\Gamma$-full label $S$ \emph{complete} whenever for each formula $A$ we have either $A \in S$ or $\neg A \in S$.
\end{definition}

Whereas for each $\Gamma \prec \Delta$ we can find a full and maximal label, not all pairs of MCSs admit a complete label.

\begin{lemma}
For each logic \ilx-that admits Figure \ref{fig:no_maximum_ex} as a model, there are MCSs $\Gamma, \Delta$ so that for no complete $S$ we have $\Gamma \sassuring{S} \Delta$.
\end{lemma}

\begin{proof}
By the assumptions on $\ilx$ the model from Figure \ref{fig:no_maximum_ex} is admissible for \ilx. Thus, we can take the corresponding \ilx-MCSs $\Gamma$ and $\Delta$ as in the proof of Lemma \ref{theorem:labelsNotClosedUnderUnions} so that $\Gamma$ is the modal theory of $w$ and $\Delta$ the modal theory of $u_1$. Let $S$ be any label with $\Gamma \sassuring{S}\Delta$. Clearly, $\neg p \notin S$ and likewise $\neg q \notin S$ since $p,q \in \Delta$. If $S$ were to be a complete label, we would need to have that $p,q \in S$ but we have seen that this is impossible.
\end{proof}

It is interesting to observe that complete labels are also MCSs so that labelled successors then become a ternary relation on MCSs. 

It seems a natural question if we can characterise which sets $\Gamma$ or which pairs of $\Gamma$ and $\Delta$ allow a complete label between them (and for which ambient logics). Another related question is for which logics there are MCSs with some particular complete label $S$ as their label.

\begin{lemma}
\label{lemma:complete_labels_form}
Let $\mathsf{Prop}$ denote the set of all propositional variables, and $\mathsf{Form}$ the set of all modal formulas.
Given any logic, any complete label equals $\widetilde{S}$ where
\[
    S = Q \cup \neg(\mathsf{Prop} \setminus Q) \cup \Box\mathsf{Form} \cup (\mathsf{Form}\rhd\mathsf{Form})
\]
for some $Q \subseteq \mathsf{Prop}$. 

\end{lemma}

\begin{proof}
    Every modal formula is a Boolean combination of propositional variables, formulas from the set  $\Box\mathsf{Form}$ and formulas from the set  $\mathsf{Form}\rhd\mathsf{Form}$, so when constructing a complete label it suffices to decide whether we will take $A$ or $\neg A$ for every formula $A$ in these three sets. A label cannot contain any formula from the set $\neg\Box\mathsf{Form}$ or the set $\neg(\mathsf{Form}\rhd\mathsf{Form})$ (see Lemma \ref{lemma:assuringnessAssures}). So, complete labels contain both $\Box\mathsf{Form}$ and $\mathsf{Form}\rhd\mathsf{Form}$. Thus, the only degree of freedom is in choosing the set $Q$ of propositional variables.
\end{proof}

The preceding lemma establishes what a complete label must look like, but does not tell us if such a label exists, i.e.~whether there is an ambient logic \ilx and \ilx-MCSs with such a label between them. We answer this in Lemma \ref{lemma:complete_labels_exist}.

\section{Existence of assuring successors}

So far, we have only been concerned with the nature and closure properties of the respecitve kind of labels. In this section we see how to guarantee the existence of informative labels.

When $\Gamma \sassuring{S} \Delta$, this enforces many formulas of the form $\neg (A \rhd B)$ to be in $\Gamma$ as we can see in the next lemma.

\begin{lemma}\label{theorem:SassuringImpliesManyNegatedRhds}
For any logic, let $\Gamma \sassuring{S} \Delta$ with $B \in \Delta$ and $S'\subseteq_{\sf fin} S$. We then have 
\[
\neg (B\rhd \bigvee_{S_i \in S'}\neg S_i) \in \Gamma.
\]
\end{lemma}

\begin{proof}
Suppose $B\rhd \bigvee_{S_i \in S'}\neg S_i \in \Gamma$. Then $\neg \neg B\rhd \bigvee_{S_i \in S'}\neg S_i \in \Gamma$, and by $\Gamma \sassuring{S} \Delta$ we would have $\neg B \in \Delta$ which is a contradiction. Thus $B\rhd \bigvee_{S_i \in S'}\neg S_i \notin \Gamma$ and by maximality 
$\neg (B\rhd \bigvee_{S_i \in S'}\neg S_i) \in \Gamma$.
\end{proof}

Conversely, the next theorem will show that given a label $S$ and maximal consistent set $\Gamma$ we have: if there are sufficiently many negated interpretability formulas related to $S$ in $\Gamma$, then we can conclude that there exists some MCS $\Delta$ with $\Gamma \sassuring{S}\Delta$.

\begin{theorem}\label{theo:ex} 
Let $\Gamma$ be an~\ilx-\mcs, $B$ a formula and $S$ a~set of formulas. If 
for any choice of finite $S' \subseteq S$ we have that  
$\neg (B\rhd \bigvee_{S_i \in S'} \neg S_i)\in\Gamma$, 
then\footnote{Lemma \ref{theorem:SassuringImpliesManyNegatedRhds} tells us that we actually have iff.} 
there exists an~\ilx-\mcs $\Delta$ such that $\Gamma\sassuring{S}\Delta\ni B, \Box\neg B$.
\end{theorem}

\begin{proof}
Suppose for a~contradiction there is no such $\Delta$. 
Then there is a~formula\footnote{There are 
finitely many $A_j$ such that there exist some formulas $S_i^j   \in S$ with 
$(\neg A_j\rhd \bigvee \neg S^j_i) \in \Gamma$ and
$\Box\neg B, B, A_j, \Box A_j \vdash \bot$. We can take $A$ to be 
$\bigwedge_j A_j$.} 
$\neg A$ such that for some $S_i \in S$ we have 
$(\neg A\rhd \bigvee \neg S_i)\in\Gamma$ and 
$B, \Box\neg B, A, \Box A \vdash \bot$. 
Then $\vdash \Box\neg B\wedge B\rhd \neg A\vee\Diamond\neg A$ and we 
get $\vdash B\rhd \neg A$. As $(\neg A\rhd \bigvee \neg S_i)\in\Gamma$, 
also $(B\rhd\bigvee \neg S_i)\in\Gamma$. A contradiction.
\end{proof}

In various applications of this theorem in modal completeness proofs we only need the lemma for a singleton set which is a direct corollary.

\begin{lemma}\label{lemm:probs}
Let $\Gamma$ be an~\ilx-\mcs such that $\neg(B\rhd C)\in\Gamma$.
Then there is an~\ilx-\mcs $\Delta$ such that $\Gamma\sassuring{ \{ \neg C \}} \Delta $ 
and $ B,\Box\neg B\in \Delta$.
\end{lemma}

\begin{proof}
Taking $S=\{ \neg C\}$ in Theorem \ref{theo:ex}.
\end{proof}

The next lemma is often given an independent proof. We show here that it in fact follows from the main theorem, Theorem \ref{theo:ex}.

\begin{lemma}\label{lemm:defies}
Let $\Gamma$ and $\Delta$ be \ilx-MCSs such that 
$A\rhd B \in \Gamma\sassuring{S}\Delta\ni A$. 
Then there is an \ilx-\mcs $\Delta'$ such that 
$\Gamma \sassuring{S} \Delta' \ni B, \Box\neg B$.
\end{lemma}

\begin{proof} First we see that for any choice of $S_i$,  
$\neg (B\rhd \bigvee \neg S_i)\in\Gamma$. Suppose not. Then for some $S_i$, $(B\rhd \bigvee \neg S_i)\in\Gamma$ because $\Gamma$ is a~MCS. 
But then $(A\rhd\bigvee \neg S_i)\in\Gamma$ and by $\Gamma\sassuring{S}\Delta$ 
we have $\neg A\in\Delta$ ($\Gamma$ is a MCS, so $(\neg \neg A\rhd\bigvee \neg S_i)\in\Gamma$). A contradiction. 
So $\neg (B\rhd \bigvee \neg S_i)\in\Gamma$ for any choice of $S_i$ and we can apply Theorem \ref{theo:ex}.
\end{proof}

The next lemma provides a sufficient condition for when an \ilx-theory can feature as label between two \ilx-MCSs.

\begin{lemma}\label{theorem:sufficientConditionForATheoryToBeLabel}
For a given logic \ilx, let $T$ be an \ilx-theory such that $T \nvdash_{\ilx} \Box \bot$. Then there are \ilx-MCSs $\Gamma$ and $\Delta$ so that $\Gamma \sassuring{T} \Delta$. 
\end{lemma}

\begin{proof}
By assumption $T \nvdash_{\ilx} \Box \bot$, so let $\Gamma$ be any \ilx-MCS that contains $T \cup \{\Diamond \top\}$.  Since $T$ is a theory, it is closed under necessitation so that $\Box T \subseteq T$, whence $\Box T \subseteq \Gamma$. Now, we consider $\neg A\rhd \neg C \in \Gamma$ with $C\in T$ (since $T$ is a theory it is closed under weakening and we may indeed take a single $C$ w.l.o.g.). But, as $\Box C \in \Box T \subseteq \Gamma$, we have $\Box A \in \Gamma$ since $\il \vdash (\neg A\rhd \neg C) \wedge \Box C \to \Box A$. Thus, any successor $\Delta$ of $\Gamma$ will contain both $A$ and $\Box A$ so for any successor $\Delta$ of $\Gamma$ we have $\Gamma \sassuring{T} \Delta$. Since $\Diamond \top \in \Gamma$ there will be indeed successors of $\Gamma$ (for example via Theorem \ref{theo:ex} taking $S=\varnothing$ and $B=\top$). 
\end{proof}

It is natural to ask if we can find a necessary and sufficient condition for some theory $T$ to feature as a label. Furthermore, in the light of Lemma \ref{theorem:informationAccrualOverFullLabels} it makes sense to ask for which chain of non-decreasing \ilx-theories $T_0\subseteq T_1 \subseteq \ldots \subseteq T_n$ we can find \ilx-MCSs $\Gamma_i$ so that
\[
\Gamma_0 \sassuring{T_0} \Gamma_1 \sassuring{T_1}  \ldots \sassuring{T_n}\Gamma_{n+1}.
\]

We now briefly return to the discussion at the end of the last section: when can a complete label appear between two MCSs. We established in Lemma \ref{lemma:complete_labels_form} the form of any such label, but does any such label exist? 

\begin{lemma}
    \label{lemma:complete_labels_exist}
    For \extil{} and some of its extensions\footnote{Some particular logics that this lemma applies to are \extil{M}, \extil{P}, \extil{W}, \extil{M_0}, and \extil{R}.} any set $S$ of the form described in Lemma \ref{lemma:complete_labels_form} is a complete label for some MCSs.
\end{lemma}
\begin{proof}
    Let $\Delta$ be the modal theory of a terminal world in a two-point regular Veltman model where this terminal world satisfies exactly the propositional variables in $S$, and let $\Gamma$ be the modal theory of its $R$-predecessor.
    It is easy to verify that $\Gamma \sassuring{S} \Delta$.
\end{proof}

The preceding lemma does not apply to, for example, any logic obtained by extending \extil{} with any single propositional variable.
In general, no set $S$ of the form described in Lemma \ref{lemma:complete_labels_form} is a complete label for all possible ambient logics, even if we restrict the claim to ambient logics that do not prove $\Box \bot$. Given such a set $S$ and the corresponding set of propositional variables $Q$, let $\il_Q$ be the smallest extension of \il with $\neg Q \cup \mathsf{Prop} \setminus Q$. The logic $\il_Q$ is a consistent logic (we can construct a Veltman model for $\il_Q$, in particular also two-point models), but no $\il_Q$-\mcs{} will contain $Q$, and thus $S$ cannot be a label between any $\il_Q$-\mcs{}s. In very general terms, such a set $S$ will be an \ilx-label if and only if for any finite collection of formulas $S_i \in S$ we have that $\extil{X} \nvdash \top \rhd \bigvee \neg S_i$ (a simple consequence of Theorem \ref{theo:ex}). 

\section{Frame conditions and labelling lemmata}\label{sec:labellinglemmata}


Although we do not yet prove completeness of any extension of \il, in this section we recall what steps there are along the way when constructing a counter-model to an unprovable formula. We can think of the step-by-step method of constructing a counter-model used e.g.\ in \cite{jogo08} for now. Later, in Section \ref{sec:ilw}, we give a proof of completeness and the finite model property for \extil{W} and a proof of completeness for \extil{P} in Section \ref{sec:seq_gen}, defining a model all at once. 

The idea in all the cases is to build a model from MCSs and define the $R$ and $S_\Gamma$ accessibility relations on them, where in particular the $R$ relation is to be defined using $\prec$. We wish to use the labels along $\prec$ to keep track of the promises posed on later added worlds by already contained interpretability formulas, and, as we shall see, also to be able to ensure we can ``locally satisfy'' the frame conditions corresponding to the additional axiom schemas, 
i.e.\ we can close the model under the characteristic property of the schema.

Let $W$ be a multiset of MCSs used in the model we wish to define. 
The main points one has to address are the following three:

\begin{enumerate}
\item\label{i:sbs:p}
  For each $\Gamma\in W$ with $\neg(A\rhd B)\in\Gamma$ we need to include a $\{\neg B\}$-assuring successor
  $\Delta$ in $W$ for which $A\in\Delta$.
\item\label{i:sbs:d}
  For each $\Gamma,\Delta\in W$ with $C\rhd D\in\Gamma \prec\Delta\ni C$ we need to include
  a $\Delta'$ in $W$ for which $\Gamma\prec\Delta'\ni D$.
  Moreover if $\Delta$ is a $T$-assuring successor of $\Gamma$
  then we should be able to choose $\Delta'$ a $T$-assuring successor of $\Gamma$
  as well (to carry promises along the $S_\Gamma$ relation).
\item\label{i:sbs:f} We need to make sure all the appropriate frame conditions are satisfied.
\end{enumerate}

The existence Lemmata \ref{lemm:probs}, \ref{lemm:defies} of the previous section ensure existence of MCSs required to witness modal formulas as specified in Item \ref{i:sbs:p} and Item \ref{i:sbs:d}. 
When working in $\il$ alone, making sure that the frame conditions are satisfied does not pose any problems \cite{joja98}, as they are just the basic properties of $R$ and $S_\Gamma$, but
with various extensions of $\il$ the situation regarding the frame conditions for the additional modal principles becomes more complicated (cf. \cite{jovelt90,jogo08}). 

\paragraph{Principle \principle{P}.}
Let us see how frame conditions locally impose requirements on MCSs, taking \extil{P} as the first example. 
The frame condition for \principle{P} is the following \cite{jovelt90}:
\[
wRw'RuS_w v \Rightarrow uS_{w'}v.
\]

The frame condition for \principle{P} imposes on MCSs the following:
\[
\Gamma \prec \Lambda \prec \Delta S_\Gamma \Sigma \Rightarrow \Delta S_{\Lambda}\Sigma.
\]

When \mcs $\Sigma\ni D$ is chosen witnessing a formula $C\rhd D\in\Gamma\prec\Delta\ni C$ by Item \ref{i:sbs:d}, we want to be able to do so in a way where not only $\Gamma\prec\Sigma$ (and the same formulas are assured), but also $\Lambda\prec\Sigma$. Moreover, if $\Lambda\prec_T\Delta$, it should be possible to choose $\Sigma$ so that $\Lambda\prec_T\Sigma$. Only then it is consistent to draw the $\Delta S_\Lambda \Sigma$ arrow required by the frame condition, as depicted in Figure \ref{fig:framecondP}. 

To see such requirements are indeed possible to meet, we will prove, for each principle, a \emph{labelling lemma}. Labelling lemmata tell us how to label the $\prec$ relation in a sufficient way to ensure we can meet the requirements imposed by frame conditions locally. 


\begin{figure}
    \centering
    \begin{minipage}[b]{0.47\textwidth} 
        \centering
        \scalebox{0.8}{
\tikzset{every picture/.style={line width=0.75pt}}
\begin{tikzpicture}[x=0.75pt,y=0.75pt,yscale=-1,xscale=1]


\draw    (63.5,779) -- (63.5,710.33) ;
\draw [shift={(63.5,708.33)}, rotate = 450] [color={rgb, 255:red, 0; green, 0; blue, 0 }  ][line width=0.75]    (10.93,-3.29) .. controls (6.95,-1.4) and (3.31,-0.3) .. (0,0) .. controls (3.31,0.3) and (6.95,1.4) .. (10.93,3.29)   ;
\draw    (63.5,708.33) -- (63.5,644) ;
\draw [shift={(63.5,642)}, rotate = 450] [color={rgb, 255:red, 0; green, 0; blue, 0 }  ][line width=0.75]    (10.93,-3.29) .. controls (6.95,-1.4) and (3.31,-0.3) .. (0,0) .. controls (3.31,0.3) and (6.95,1.4) .. (10.93,3.29)   ;
\draw    (63.5,642) .. controls (99.95,622.63) and (127.66,666) .. (156.2,644.05) ;
\draw [shift={(157.5,643)}, rotate = 500] [color={rgb, 255:red, 0; green, 0; blue, 0 }  ][line width=0.75]    (10.93,-3.29) .. controls (6.95,-1.4) and (3.31,-0.3) .. (0,0) .. controls (3.31,0.3) and (6.95,1.4) .. (10.93,3.29)   ;
\draw    (63.5,779) -- (156.36,644.65) ;
\draw [shift={(157.5,643)}, rotate = 484.65] [color={rgb, 255:red, 0; green, 0; blue, 0 }  ][line width=0.75]    (10.93,-3.29) .. controls (6.95,-1.4) and (3.31,-0.3) .. (0,0) .. controls (3.31,0.3) and (6.95,1.4) .. (10.93,3.29)   ;
\draw  [dash pattern={on 4.5pt off 4.5pt}]  (63.5,642) .. controls (110.03,598.77) and (99.71,636.56) .. (155.78,642.82) ;
\draw [shift={(157.5,643)}, rotate = 185.58] [color={rgb, 255:red, 0; green, 0; blue, 0 }  ][line width=0.75]    (10.93,-3.29) .. controls (6.95,-1.4) and (3.31,-0.3) .. (0,0) .. controls (3.31,0.3) and (6.95,1.4) .. (10.93,3.29)   ;
\draw  [dash pattern={on 4.5pt off 4.5pt}]  (63.5,708.33) -- (155.86,644.14) ;
\draw [shift={(157.5,643)}, rotate = 505.2] [color={rgb, 255:red, 0; green, 0; blue, 0 }  ][line width=0.75]    (10.93,-3.29) .. controls (6.95,-1.4) and (3.31,-0.3) .. (0,0) .. controls (3.31,0.3) and (6.95,1.4) .. (10.93,3.29)   ;

\draw (53,776.67) node   [align=left] {$\displaystyle \Gamma $};
\draw (53,710.67) node   [align=left] {$\displaystyle \Lambda $};
\draw (54,645.67) node   [align=left] {$\displaystyle \Delta $};
\draw (168,644.67) node   [align=left] {$\displaystyle \Sigma $};
\draw (102,649.67) node   [align=left] {$\displaystyle S_{\Gamma }$};
\draw (53.36,679) node   [align=left] {$\displaystyle T$};
\draw (52.36,749) node   [align=left] {$\displaystyle S$};
\draw (142.36,714) node   [align=left] {$\displaystyle S\cup \Lambda ^{\boxdot }_{T}$};
\draw (102,607.67) node   [align=left] {$\displaystyle S_{\Lambda }$};
\draw (91.36,678) node   [align=left] {$\displaystyle T$};

\end{tikzpicture}
}
        \caption{\footnotesize Ensuring frame condition for  ${\sf P}$}
        \label{fig:framecondP}
    \end{minipage}
    \begin{minipage}[b]{0.47\textwidth} 
        \centering
        \scalebox{0.8}{
\begin{tikzpicture}[x=0.75pt,y=0.75pt,yscale=-1,xscale=1]

\draw    (262.5,788) -- (262.5,719.33) ;
\draw [shift={(262.5,717.33)}, rotate = 450] [color={rgb, 255:red, 0; green, 0; blue, 0 }  ][line width=0.75]    (10.93,-3.29) .. controls (6.95,-1.4) and (3.31,-0.3) .. (0,0) .. controls (3.31,0.3) and (6.95,1.4) .. (10.93,3.29)   ;
\draw    (356.5,718.33) -- (356.5,654) ;
\draw [shift={(356.5,652)}, rotate = 450] [color={rgb, 255:red, 0; green, 0; blue, 0 }  ][line width=0.75]    (10.93,-3.29) .. controls (6.95,-1.4) and (3.31,-0.3) .. (0,0) .. controls (3.31,0.3) and (6.95,1.4) .. (10.93,3.29)   ;
\draw    (262.5,717.33) .. controls (298.95,697.96) and (326.66,741.33) .. (355.2,719.38) ;
\draw [shift={(356.5,718.33)}, rotate = 500] [color={rgb, 255:red, 0; green, 0; blue, 0 }  ][line width=0.75]    (10.93,-3.29) .. controls (6.95,-1.4) and (3.31,-0.3) .. (0,0) .. controls (3.31,0.3) and (6.95,1.4) .. (10.93,3.29)   ;
\draw    (262.5,788) -- (354.89,719.52) ;
\draw [shift={(356.5,718.33)}, rotate = 503.46] [color={rgb, 255:red, 0; green, 0; blue, 0 }  ][line width=0.75]    (10.93,-3.29) .. controls (6.95,-1.4) and (3.31,-0.3) .. (0,0) .. controls (3.31,0.3) and (6.95,1.4) .. (10.93,3.29)   ;
\draw  [dash pattern={on 4.5pt off 4.5pt}]  (262.5,717.33) -- (354.86,653.14) ;
\draw [shift={(356.5,652)}, rotate = 505.2] [color={rgb, 255:red, 0; green, 0; blue, 0 }  ][line width=0.75]    (10.93,-3.29) .. controls (6.95,-1.4) and (3.31,-0.3) .. (0,0) .. controls (3.31,0.3) and (6.95,1.4) .. (10.93,3.29)   ;

\draw (252,785.67) node   [align=left] {$\displaystyle \Gamma $};
\draw (252,719.67) node   [align=left] {$\displaystyle \Delta $};
\draw (367,653.67) node   [align=left] {$\displaystyle \Omega $};
\draw (301,725.67) node   [align=left] {$\displaystyle S_{\Gamma }$};
\draw (234.36,750) node   [align=left] {$\displaystyle S\cup \Delta ^{\square }$};
\draw (335.36,761) node   [align=left] {$\displaystyle S\cup \Delta ^{\square }$};
\draw (369,718.67) node   [align=left] {$\displaystyle \Sigma $};

\end{tikzpicture}
}
        \caption{\footnotesize Ensuring frame condition for  ${\sf M}$}
        \label{fig:framecondM}
    \end{minipage}
\end{figure}   

\begin{lemma}\label{lemm:labp}
For logics containing ${\sf P}$ we have
$\Gamma\sassuring{S}\Lambda\sassuring{T}\Delta\Rightarrow\Gamma\sassuring{S\cup\boxdotset{T}{\Lambda}}\Delta$.
\end{lemma}

\begin{proof}
Suppose $C\rhd\bigvee\neg S_i\vee\bigvee A_j\vee\Diamond A_j\in\Gamma$,
where $\Box\neg A_j,\neg A_j\in\boxdotset{T}{\Lambda}$.
Then $C\rhd\bigvee\neg S_i\vee\bigvee A_j\in\Gamma$ and thus
by ${\sf P}$ we obtain $C\rhd\bigvee\neg S_i\vee\bigvee A_j\in\Lambda$.
Since $\Gamma\sassuring{S}\Lambda$ we have $\Box\bigwedge S_i\in\Lambda$ so we 
obtain
$C\rhd\bigvee A_j\in\Lambda$. But for each $A_j$ we have 
$A_j\rhd\bigvee\neg T_{jk}\in\Lambda$
and thus $C\rhd\bigvee\neg T_{jk}\in\Lambda$. Since $\Lambda\prec_T\Delta$ we 
conclude
$\neg C,\Box\neg C\in\Delta$.
\end{proof}


In the case of P, a simpler labelling lemma can be used to ensure the frame condition locally, provided we consider the labels that are \emph{full} ($S$ a $\Gamma$-full label, and $T$ a $\Lambda$-full label).

\begin{lemma}\label{lemm:labfullp}
For logics containing ${\sf P}$ we have 
$$
\Gamma\sassuring{S}\Lambda\sassuring{T}\Delta\Rightarrow\Gamma\sassuring{T}\Delta
$$
\end{lemma}

\begin{proof}
Assume $\Gamma\sassuring{S}\Lambda\sassuring{T}\Delta$, and $C\rhd \bigvee \neg T_i \in\Gamma$. Then by ${\sf P}$ we know $C\rhd \bigvee \neg T_i \in\Lambda$. Since $\Lambda\prec_T\Delta$ we conclude $\neg C,\Box\neg C\in\Delta$. 
\end{proof}

Note that the lemma is true in the case of ordinary labels, but in that case, the previous lemma gives us more precise labelling information to ensure the frame condition locally. This is because only for full labels we in fact have $S\cup\boxdotset{T}{\Lambda}\subseteq T$.  

\paragraph{Principle \principle{M}.}
The frame condition for \principle{M} is the following \cite{jovelt90}:
\[
wRuS_w v Rz \Rightarrow uRz.
\]

The frame condition for \principle{M} imposes on MCSs the following:
\[
\Delta S_\Gamma \Sigma \prec \Omega \Rightarrow \Delta \prec\Omega.
\]

When \mcs $\Sigma\ni D$ is chosen witnessing a formula $C\rhd D\in\Gamma\prec\Delta\ni C$ by Item \ref{i:sbs:d}, we want to do so in such a way that whenever we later need to add a \mcs $\Omega$ with $\Sigma\prec\Omega$, we can also draw the $\Delta\prec\Omega$ arrow. Therefore we need to ensure $\boxset{\emptyset}{\Delta}$ along the $\Gamma\prec\Sigma$ arrow (as we remarked previously, one can think of the set $\boxset{\emptyset}{\Delta}$ as simply $\{ \Box C \mid \Box C \in \Delta \}$), we achieve this by ensuring $\boxset{\emptyset}{\Delta}$ along the $\Gamma\prec\Delta$ arrow. The situation is depicted in 
Figure \ref{fig:framecondM}. The corresponding labelling lemma is the following:

\begin{lemma}\label{lemm:labm}
For logics containing \principle{M} we have $\Gamma \sassuring{S} \Delta 
\Rightarrow \Gamma \sassuring{S\cup \boxset{\emptyset}{\Delta}} \Delta$.
\end{lemma}

\begin{proof}
Assume that for some $\Box C_j \in \boxset{\emptyset}{\Delta}$ we have 
$(A\rhd\bigvee \neg S_i \vee \bigvee \neg \Box C_j) \in \Gamma$.
By \principle{M},
$(A \wedge \bigwedge \Box C_j \rhd \bigvee \neg S_i) \in \Gamma$,
whence $\boxdot \neg (A \wedge \bigwedge \Box C_j) \in \Delta$.
As $\bigwedge \Box C_j \in \Delta$, we conclude
$\neg A, \Box \neg A \in \Delta$.
\end{proof}

In the case of \principle{M}, we have no simpler labelling lemma in case $S$ is a $\Gamma$-full label.

\paragraph{Principle \principle{M_0}.}
The frame condition for \principle{M_0} is the following \cite{jogo08}:
\[
wRuRxS_w v Rz \Rightarrow uRz
\]

The frame condition for \principle{M_0} imposes on \mcs the following:
\[
\Gamma \prec \Delta \prec \Delta' S_\Gamma \Sigma \prec \Omega \Rightarrow \Delta \prec \Omega.
\]

When \mcs $\Sigma\ni D$ is chosen witnessing a formula $C\rhd D\in\Gamma\prec\Delta\prec\Delta'\ni C$ by Item \ref{i:sbs:d}, we want to do so in such a way that whenever we later need to add a \mcs $\Omega$ with $\Sigma\prec\Omega$, we can also draw the $\Delta\prec\Omega$ arrow. Therefore we again need to ensure $\boxset{\emptyset}{\Delta}$ along the $\Gamma\prec\Sigma$ arrow. The situation is depicted in Figure \ref{fig:framecondM_0}, and the corresponding labelling lemma is the following (as before, we do not have a special lemma in case the labels are full): 

\begin{lemma}\label{lemm:labmo}
For logics containing \principle{M_0} we have
$\Gamma \sassuring{S} \Delta \sucs \Delta' \Rightarrow 
\Gamma \sassuring{S \cup \boxset{\emptyset}{\Delta}} \Delta'$.
\end{lemma}

\begin{proof}
Suppose $C\rhd\bigvee \neg S_i\vee\bigvee\Diamond A_j\in\Gamma$,
where $\Box\neg A_j\in\boxset{\emptyset}{\Delta}$.
By ${\sf M}_0$ we obtain $\Diamond C\wedge\bigwedge\Box\neg A_j\rhd\bigvee 
\neg S_i\in\Gamma$.
So, since $\Gamma\sassuring{S}\Delta$ and $\bigwedge\Box\neg A_j\in\Delta$ we 
obtain $\Box\neg C\in\Delta$ and thus $\Box\neg C,\neg C\in\Delta'$.
\end{proof}



\begin{figure}
    \centering
    \begin{minipage}[b]{0.47\textwidth} 
        \centering
        \scalebox{0.8}{
\begin{tikzpicture}[x=0.75pt,y=0.75pt,yscale=-1,xscale=1]

\draw    (452.5,784) -- (452.5,715.33) ;
\draw [shift={(452.5,713.33)}, rotate = 450] [color={rgb, 255:red, 0; green, 0; blue, 0 }  ][line width=0.75]    (10.93,-3.29) .. controls (6.95,-1.4) and (3.31,-0.3) .. (0,0) .. controls (3.31,0.3) and (6.95,1.4) .. (10.93,3.29)   ;
\draw    (452.5,713.33) -- (452.5,649) ;
\draw [shift={(452.5,647)}, rotate = 450] [color={rgb, 255:red, 0; green, 0; blue, 0 }  ][line width=0.75]    (10.93,-3.29) .. controls (6.95,-1.4) and (3.31,-0.3) .. (0,0) .. controls (3.31,0.3) and (6.95,1.4) .. (10.93,3.29)   ;
\draw    (452.5,647) .. controls (488.95,627.63) and (516.66,671) .. (545.2,649.05) ;
\draw [shift={(546.5,648)}, rotate = 500] [color={rgb, 255:red, 0; green, 0; blue, 0 }  ][line width=0.75]    (10.93,-3.29) .. controls (6.95,-1.4) and (3.31,-0.3) .. (0,0) .. controls (3.31,0.3) and (6.95,1.4) .. (10.93,3.29)   ;
\draw    (452.5,784) -- (545.36,649.65) ;
\draw [shift={(546.5,648)}, rotate = 484.65] [color={rgb, 255:red, 0; green, 0; blue, 0 }  ][line width=0.75]    (10.93,-3.29) .. controls (6.95,-1.4) and (3.31,-0.3) .. (0,0) .. controls (3.31,0.3) and (6.95,1.4) .. (10.93,3.29)   ;
\draw  [dash pattern={on 4.5pt off 4.5pt}]  (452.5,713.33) -- (545.34,583.29) ;
\draw [shift={(546.5,581.67)}, rotate = 485.52] [color={rgb, 255:red, 0; green, 0; blue, 0 }  ][line width=0.75]    (10.93,-3.29) .. controls (6.95,-1.4) and (3.31,-0.3) .. (0,0) .. controls (3.31,0.3) and (6.95,1.4) .. (10.93,3.29)   ;
\draw    (546.5,648) -- (546.5,583.67) ;
\draw [shift={(546.5,581.67)}, rotate = 450] [color={rgb, 255:red, 0; green, 0; blue, 0 }  ][line width=0.75]    (10.93,-3.29) .. controls (6.95,-1.4) and (3.31,-0.3) .. (0,0) .. controls (3.31,0.3) and (6.95,1.4) .. (10.93,3.29)   ;

\draw (442,781.67) node   [align=left] {$\displaystyle \Gamma $};
\draw (442,715.67) node   [align=left] {$\displaystyle \Delta $};
\draw (443,650.67) node   [align=left] {$\displaystyle \Delta '$};
\draw (557,649.67) node   [align=left] {$\displaystyle \Sigma $};
\draw (478,628.67) node   [align=left] {$\displaystyle S_{\Gamma }$};
\draw (441.36,748) node   [align=left] {$\displaystyle S$};
\draw (531.36,719) node   [align=left] {$\displaystyle S\cup \Delta ^{\square }$};
\draw (560,582.67) node   [align=left] {$\displaystyle \Omega $};

\end{tikzpicture}
}
        \caption{\footnotesize Ensuring frame condition for  ${\sf M}_0$}
        \label{fig:framecondM_0}
    \end{minipage}
    \begin{minipage}[b]{0.47\textwidth} 
        \centering
        \scalebox{0.8}{
\begin{tikzpicture}[x=0.75pt,y=0.75pt,yscale=-1,xscale=1]

\draw    (57.5,1039) -- (57.5,970.33) ;
\draw [shift={(57.5,968.33)}, rotate = 450] [color={rgb, 255:red, 0; green, 0; blue, 0 }  ][line width=0.75]    (10.93,-3.29) .. controls (6.95,-1.4) and (3.31,-0.3) .. (0,0) .. controls (3.31,0.3) and (6.95,1.4) .. (10.93,3.29)   ;
\draw    (57.5,968.33) -- (57.5,904) ;
\draw [shift={(57.5,902)}, rotate = 450] [color={rgb, 255:red, 0; green, 0; blue, 0 }  ][line width=0.75]    (10.93,-3.29) .. controls (6.95,-1.4) and (3.31,-0.3) .. (0,0) .. controls (3.31,0.3) and (6.95,1.4) .. (10.93,3.29)   ;
\draw    (57.5,902) .. controls (93.95,882.63) and (121.66,926) .. (150.2,904.05) ;
\draw [shift={(151.5,903)}, rotate = 500] [color={rgb, 255:red, 0; green, 0; blue, 0 }  ][line width=0.75]    (10.93,-3.29) .. controls (6.95,-1.4) and (3.31,-0.3) .. (0,0) .. controls (3.31,0.3) and (6.95,1.4) .. (10.93,3.29)   ;
\draw    (57.5,1039) -- (150.36,904.65) ;
\draw [shift={(151.5,903)}, rotate = 484.65] [color={rgb, 255:red, 0; green, 0; blue, 0 }  ][line width=0.75]    (10.93,-3.29) .. controls (6.95,-1.4) and (3.31,-0.3) .. (0,0) .. controls (3.31,0.3) and (6.95,1.4) .. (10.93,3.29)   ;
\draw  [dash pattern={on 4.5pt off 4.5pt}]  (57.5,968.33) -- (150.34,838.29) ;
\draw [shift={(151.5,836.67)}, rotate = 485.52] [color={rgb, 255:red, 0; green, 0; blue, 0 }  ][line width=0.75]    (10.93,-3.29) .. controls (6.95,-1.4) and (3.31,-0.3) .. (0,0) .. controls (3.31,0.3) and (6.95,1.4) .. (10.93,3.29)   ;
\draw    (151.5,903) -- (151.5,838.67) ;
\draw [shift={(151.5,836.67)}, rotate = 450] [color={rgb, 255:red, 0; green, 0; blue, 0 }  ][line width=0.75]    (10.93,-3.29) .. controls (6.95,-1.4) and (3.31,-0.3) .. (0,0) .. controls (3.31,0.3) and (6.95,1.4) .. (10.93,3.29)   ;
\draw  [dash pattern={on 4.5pt off 4.5pt}]  (57.5,902) .. controls (81.26,823.13) and (104.04,877.23) .. (150.1,837.89) ;
\draw [shift={(151.5,836.67)}, rotate = 498.44] [color={rgb, 255:red, 0; green, 0; blue, 0 }  ][line width=0.75]    (10.93,-3.29) .. controls (6.95,-1.4) and (3.31,-0.3) .. (0,0) .. controls (3.31,0.3) and (6.95,1.4) .. (10.93,3.29)   ;

\draw (47,1036.67) node   [align=left] {$\displaystyle \Gamma $};
\draw (47,970.67) node   [align=left] {$\displaystyle \Lambda $};
\draw (48,905.67) node   [align=left] {$\displaystyle \Delta $};
\draw (162,904.67) node   [align=left] {$\displaystyle \Sigma $};
\draw (91,888.67) node   [align=left] {$\displaystyle S_{\Gamma }$};
\draw (46.36,1003) node   [align=left] {$\displaystyle S$};
\draw (136.36,974) node   [align=left] {$\displaystyle S\cup \Lambda ^{\square }_{T}$};
\draw (165,837.67) node   [align=left] {$\displaystyle \Omega $};
\draw (46.36,941) node   [align=left] {$\displaystyle T$};
\draw (98.36,927) node   [align=left] {$\displaystyle T$};
\draw (96,837.67) node   [align=left] {$\displaystyle S_{\Lambda }$};

\end{tikzpicture}
}
        \caption{\footnotesize Ensuring frame condition for  ${\sf R}$}
        \label{fig:framecondR}
    \end{minipage}
\end{figure}  

\paragraph{Principle \principle{R}.}
Last we will look at a more complicated case of \extil{R}.

The frame condition for the principle \principle{R} is the following \cite{jogo11}.\footnote{
In \cite{jogo04} the modal principle
$A\rhd B\rightarrow\neg(A\rhd\neg C)\wedge(D\rhd C)\rhd B\wedge\Box C$
was called ${\sf R}$. This principle and the one called ${\sf R}$ in this
paper are easily seen to be equivalent over \il.}
\[
wRxRyS_wy'Rz\Rightarrow yS_xz.
\]

On \mcs, the condition imposes the following:
\[
\Gamma \prec \Lambda \prec \Delta S_\Gamma \Sigma \prec \Omega\Rightarrow \Delta S_\Lambda \Omega.
\]

The frame condition is depicted in Figure \ref{fig:framecondR}. 
Assume $\Sigma\ni D$ was chosen as a witness for $C\rhd D\in\Gamma R \Delta\ni C$.
Since $\Delta$ lies $T$-assuring above $\Lambda$,
we should not only make sure that $\Sigma$ lies $S$-assuring above $\Gamma$, but also that any successor $\Omega$ of $\Sigma$ lies $T$-assuring above $\Lambda$. Only then we would be justified to draw the required $\Delta S_\Lambda \Omega$ arrow. One way to guarantee $\Lambda\sassuring{T}\Omega$ is to ensure $\boxset{T}{\Lambda}$ along the $\Gamma\prec\Sigma$ arrow: whenever $B\rhd \bigvee\neg T_i\in\Lambda$, we have $\Box \neg B\in \boxset{T}{\Lambda}$ and this puts $\Box\neg B \in \Sigma$ and $\Box\neg B,\neg B \in \Omega$ as required.

The corresponding labelling lemma is the following:

\begin{lemma}\label{lemm:labeltrans}\label{lemm:labr}
For logics containing \principle{R} we have
$\Gamma \sassuring{S} \Lambda \sassuring{T} \Delta \Rightarrow 
\Gamma \sassuring{S \cup \boxset{T}{\Lambda}} \Delta$.
\end{lemma}

\begin{proof}
We consider $A$ such that for some $S_i \in S$ and 
some $\Box \neg A_j \in \boxset{T}{\Lambda}$, we have
$(A \rhd \bigvee \neg S_i \vee \bigvee \Diamond A_j) \in \Gamma$.
By \principle{R} we obtain 
$(\neg (A \rhd \bigvee A_j)\rhd \bigvee \neg S_i)\in \Gamma$,
thus by $\Gamma \sassuring{S} \Lambda$ we get
$(A \rhd \bigvee A_j) \in \Lambda$.
As $(A_j \rhd \bigvee \neg T_{kj})\in \Lambda$, also
$(A \rhd \bigvee \neg T_{kj})\in \Lambda$. By 
$\Lambda \sassuring{T} \Delta$ we conclude $\boxdot \neg A \in \Delta$.
\end{proof}

In the case of \principle{R}, a simpler labelling lemma can be used to ensure the frame condition locally if $T$ is $\Lambda$-full: 

\begin{lemma}\label{lemm:labfullr}
For logics containing \principle{R} we have 
$\Gamma \sassuring{S} \Lambda \sassuring{T} \Delta \Rightarrow\Gamma\sassuring{S \cup \Box T} \Delta$
\end{lemma}

\begin{proof}
Assume $A\rhd \bigvee \neg S_i \vee \bigvee \neg\Box T_j \in \Gamma$. Then, by \principle{R}, we obtain $\neg(A\rhd \bigvee \neg T_j)\rhd \bigvee\neg S_i \in \Gamma$ and by $\Gamma\sassuring{S}\Lambda$ we know $\boxdot(A\rhd \bigvee \neg T_j)\in \Lambda$, and $\boxdot \neg A\in \Delta$ as required.
\end{proof}

As before in the case of logics containing \principle{P} and Lemma \ref{lemm:labfullp}, this lemma ensures the frame condition locally provided the labels are full: for in this case $\boxdotset{T}{\Lambda}\subseteq T$ and therefore, because $T$ is a theory, $\boxset{T}{\Lambda}\subseteq \Box T$, and consequently $S\cup\boxset{T}{\Lambda}\subseteq S\cup\Box T$.
Thus sufficient information is carried by the composed label.


\paragraph{Case of \extil{W}.}
For this logic we did not manage to find a labelling lemma.
Instead, let us state two existence lemmata for \extil{W}, a logic for which only second order frame properties are known (\cite{jogo11}, \cite{Mikec-Perkov-Vukovic-17}).

\begin{lemma}\label{lemm:wprobs}
Suppose $\neg(A\rhd B)\in\Gamma$ and $\Gamma$ is an \ilx-\mcs. There exists some \ilx-\mcs $\Delta$ with
$\Gamma\prec_{\{\Box\neg A,\neg B\}}\Delta$ and $A\in\Delta$.
\end{lemma}

\begin{proof}
Suppose for a contradiction that there is no such $\Delta$.
Then there are finitely many formulas $E_i$ such that $(E_i\rhd \Diamond A\vee B)\in\Gamma$ and $A, \{\neg E_i, \Box\neg E_i\}_i \vdash\bot$. Let $E = \bigvee_i E_i$. By \extil{} and maximal consistency we have $(E\rhd \Diamond A\vee B)\in\Gamma$ and
$A, \neg E, \Box\neg E\vdash\bot$. Thus $\vdash A\rhd E$.
Then $(A\rhd \Diamond A\vee B)\in\Gamma$ and by the principle $\principle{W}$ we have $A\rhd B\in\Gamma $.
The contradiction.
\end{proof}

\begin{lemma}\label{lemm:wdefies}
For logics containing \principle{W} we have that if
$B\rhd C \in \Gamma \sassuring{S} \Lambda \ni B$ then there exists
$\Delta$ (an \ilx-\mcs w.r.t.~the same logic \ilx) with
$\Gamma \sassuring{S\cup \{ \Box \neg B\}} \Delta \ni C, \Box \neg C$.
\end{lemma}

\begin{proof}
Suppose for a contradiction that no such $\Delta$ exists. Then
for some formula $A$ with 
$(A \rhd \bigvee \neg S_i \vee \Diamond B) \in \Gamma$, we get
$C , \Box \neg C, \neg A, \Box \neg A \vdash \bot$, whence
$\vdash C\rhd A$. Thus 
$B\rhd C \rhd A \rhd \bigvee \neg S_i \vee \Diamond B \in \Gamma$.
By \principle{W}, $B \rhd \bigvee \neg S_i \in \Gamma$ which 
contradicts 
$\Gamma\sassuring{S} \Lambda \ni B$.
\end{proof}


\section{Going finite}\label{sec:finite}
Proving the decidability of an interpretability logic is 
in all known cases done by either proving completeness with respect to finite models (\cite{jovelt90}, \cite{jonvelt99}) or by showing (in case completeness is known w.r.t.\ some class of models) that for any satisfiable formula there is a finite model satisfying it  (\cite{Perkov-Vukovic-16}, \cite{Mikec-Perkov-Vukovic-17}). 
The finite model property is easier to achieve 
if the building blocks of the model are finite sets instead of
infinite maximal consistent sets.

In completeness proofs, a turn that is  usually made to obtain finite building blocks is to
work with truncated parts of maximal consistent sets. These truncated parts should
be large enough to allow for the basic reasoning, and this gives 
rise to the notion of so-called adequate sets. Note that different 
logics yield different notions of adequacy. In order to obtain
the finite model property along with modal completeness of \extil{W},
in the next section we will use the following notion of adequacy.

\begin{definition}[Adequate set]
We say that a set of formulas $\Phi$ is \emph{adequate} iff
\begin{enumerate}

\item
  $\bot\rhd\bot\in\Phi$;

\item
  $\Phi$ is closed under single negation and subformulas;

\item
  If both $A$ is an antecedent or consequent of some $\rhd$ formula
  in $\Phi$ and so is $B$ then $A\rhd B\in\Phi$.

\end{enumerate}
\end{definition}

It is clear that any formula is contained in some finite and minimal adequate set.
For a formula $F$ we will denote this set by $\Phi(F)$.
Here and in the following section MCSs are subsets of, and maximal w.r.t., some adequate set $\Phi$.
Since our maximal consistent sets are more restricted we should also modify the notion
of an assuring successor a bit.

\newcommand{\fsassuring}[1]{\prec_{#1}^\Phi}

\begin{definition}[\pair{S,\Phi}-assuring successor]
Let $\Phi$ be a finite adequate set, $S\subseteq\Phi$ and $\Gamma,\Delta\subseteq\Phi$
be maximal consistent sets. We say that $\Delta$ is an $\pair{S,\Phi}$-assuring successor
of $\Gamma$ ($\Gamma\fsassuring{S}\Delta$) iff for each $\Box\neg A\in\Phi$ we have
\[
\Gamma\vdash A\rhd\bigvee_{S_i\in S}\neg S_i\Rightarrow \neg A,\Box\neg A\in\Delta,
\]
and if moreover for some $\Box C\in\Delta$ we have $\Box C\not\in\Gamma$.
\end{definition}

Note that by the requirement $\Box\neg A\in\Phi$ the usual reading of $\prec$ in
extensions of \gl coincides with $\fsassuring{\emptyset}$.
So we will write $\prec$ for $\fsassuring{\emptyset}$.
The following two lemmas follow from their infinite counterparts (Lemma \ref{lemm:wprobs} and Lemma \ref{lemm:wdefies}), by taking intersections between the sets given by those lemmas, and the set $\Phi$.

\begin{lemma}\label{lemm:wfin:probs}
Let $\Gamma\subseteq\Phi$ be maximal consistent.
If $\neg(A\rhd B)\in\Gamma$ then there exists some maximal consistent set $\Delta\subseteq\Phi$
such that $A\in\Delta$ and $\Gamma\fsassuring{\{\neg B,\Box\neg A\}}\Delta$.
\end{lemma}
\begin{proof}
    Since $\Gamma$ is consistent, there is an extension $\Gamma' \supseteq \Gamma$ that is maximal consistent w.r.t.\ the set of all modal formulas (not just $\Phi$).
    By Lemma \ref{lemm:wprobs}, there is a set $\Delta'$ with $\Gamma' \prec_{ \{\Box \neg A, \neg B\} } \Delta' \ni A$, and $\Delta'$ is maximal consistent w.r.t.\ the set of all modal formulas.
    The set $\Delta = \Delta' \cap \Phi$ is the required MCS.
\end{proof}


\begin{lemma}\label{lemm:wfin:defies}
Let $\Gamma,\Delta\subseteq\Phi$ be maximal consistent and $S\subseteq\Phi$.
If $A\rhd B\in\Gamma$, $\Gamma\fsassuring{S}\Delta$ and $A\in\Delta$ then there exists
some maximal consistent $\Delta'\subseteq\Phi$ with $B\in\Delta'$ and
$\Gamma\fsassuring{S\cup\{\Box\neg A\}}\Delta'$.
\end{lemma}
\begin{proof}
    Analogous to the proof of Lemma \ref{lemm:wfin:probs}, but this time employing Lemma \ref{lemm:wdefies}.
\end{proof}

\section{The logic \ilw}\label{sec:ilw}
\newcommand{\tuple}[1]{\langle #1 \rangle}

As a demonstration of the use of assuringness we will give in this section
a relatively simple proof of the known fact that \extil{W} is a complete logic.

In what follows we let $\Phi$ be some fixed finite adequate set and reason with $\extil{W}$
(e.g. $\vdash$ is $\extil{W}$-provable, and consistent is $\extil{W}$-consistent).
The rest of this section is devoted to the proof of the following theorem.

\begin{theorem}[Completeness of \extil{W} \cite{jonvelt99}]
\label{thm:ilw-complete}
\extil{W} is complete with respect to finite Veltman frames $\langle W,R,S\rangle$
in which, for each $w\in W$, the relation $R \circ S_w$ is conversely well-founded.
\end{theorem}

Suppose $\not\vdash G$. Let $\Phi=\Phi(\neg G)$
and let $\Gamma\subseteq\Phi$ be a maximal consistent set that contains $\neg G$.
We will construct a Veltman model $\langle W,R,\{ S_w : w \in W \},V\rangle$
in which for each $w\in W$ we have that $R \circ S_w$ is conversely well-founded.
Each $w\in W$ will be a tuple the second component of which---denoted by $(w)_1$---will be a maximal consistent subset of $\Phi$.
For some $w\in W$ we will have $(w)_1=\Gamma$
and we will finish the proof by proving a
truth lemma: $w\Vdash A$ iff $A\in(w)_1$.

Let the \emph{height} of a maximal consistent $\Delta\subseteq\Phi$ be defined as the
number of $\Box$-formulas in $\Delta$ minus the number of $\Box$-formulas in $\Gamma$.
For sequences $\sigma_0$ and $\sigma_1$ we write $\sigma_0\subseteq\sigma_1$ iff
$\sigma_0$ is an initial, but not necessarily proper subsequence of $\sigma_1$.
For two sequences $\sigma_0$ and $\sigma_1$, $\sigma_0 * \sigma_1$ denotes the concatenation
of the two sequences. If $S$ is a set of formulas then $\tuple{S}$ is the sequence
of length one and only element $S$.
Let us now define $\langle W,R,\{ S_w : w \in W \},V\rangle$.

\begin{enumerate}

\item $W$ is the set of pairs $\tuple{\sigma,\Delta}$ where
  $\Delta\subseteq\Phi$ is maximal consistent such that either $\Gamma=\Delta$ or $\Gamma\prec\Delta$ and
  $\sigma$ is a finite sequence of subsets of $\Phi$ the length of which does not exceed the height of
  $\Delta$.
  For $w=\tuple{\sigma,\Delta}$ we write $(w)_0$ for $\sigma$ and $(w)_1$ for $\Delta$.

\item $wRv$ iff 
  for some $S$ we have $(v)_0\supseteq (w)_0*\tuple{S}$ and $(w)_1\fsassuring{S}(v)_1$.

\item $xS_wy$ iff $wRx,y$ and, $xRy$ or $x=y$ or both \ref{i:sdef:2} and \ref{i:sdef:3} hold:
  \begin{enumerate}
    \item\label{i:sdef:2} If $(x)_0=(w)_0 * \tuple{S} * \tau_x$ and
      $(y)_0=(w)_0 * \tuple{T} * \tau_y$, then $S\subseteq T$.
    \item\label{i:sdef:3}
      For some $C\rhd D\in (w)_1$ we have $\Box\neg C\in T$ and, $C\in(x)_1$ or $\Diamond C\in(x)_1$.
  \end{enumerate}

\item $V(p)=\{ w\in W\mid p\in (w)_1 \}$.

\end{enumerate}

We shall now see that this defines an \extil{W}-model. First we see that $R$ behaves properly.

\begin{lemma}
$R$ is transitive and conversely well-founded.
\end{lemma}

\begin{proof}
Transitivity follows from the fact that $(x)_1\fsassuring{S}(y)_1\prec(z)_1$
implies $(x)_1\fsassuring{S}(z)_1$.
Converse well-foundedness now follows from the fact that our model is finite and $R$ is irreflexive.
\end{proof}

Next we show that the $S_w$-relations comply with their requirements.

\begin{lemma}
$wRxRy$ implies $xS_wy$.
Also $wRx$ implies $xS_wx$.
Finally, $S_w$ is transitive.
\end{lemma}

\begin{proof}

The first two assertions hold by definition.
So suppose $xS_wyS_wz$. Let us fix
$(x)_0 \supseteq (w)_0 * \tuple{S}$,
$(y)_0 \supseteq (w)_0 * \tuple{T}$ and
$(z)_0 \supseteq (w)_0 * \tuple{U}$.
We distinguish two cases.

Case 1: $xRy$ or $x=y$. If $x=y$ then we are done so we assume $xRy$.
If $yRz$ or $y=z$ then we are also easily done. So, we assume that
for some $C\rhd D\in(w)_1$ we have $\Box\neg C\in U$ and,
$C\in(y)_1$ or $\Diamond C\in(y)_1$. Since $(x)_1\prec (y)_1$ we have
that $\Diamond C \in(x)_1$ and thus we conclude $xS_wz$.

Case 2: $\neg xRy$ and $x\neq y$. In this case there exists some
$C\rhd D\in(w)_1$ with $\Box\neg C\in T$ and $C\in(x)_1$ or $\Diamond C\in(x)_1$.
Whatever the reason for $yS_wz$ is, we always have $T\subseteq U$ and thus $\Box\neg C\in U$.
So we conclude $xS_wz$.
\end{proof}

Finally, we check the frame condition for \textsf{W}: that the relation $R \circ S_w$ is conversely well-founded

\begin{lemma}
The relation $R \circ S_w$ is conversely well-founded.
\end{lemma}
\begin{proof}

Suppose we have an infinite sequence
\[
x_0S_wy_0Rx_1S_wy_1R\cdots.
\]
For each $i\geq 0$, fix $X_i$ and $Y_i$ such that
$(x_i)_0\supseteq(w)_0*\tuple{X_i}$ and
$(y_i)_0\supseteq(w)_0*\tuple{Y_i}$.
We may assume that, for each $i$, $x_i\neq y_i$ and $\neg x_iRy_i$.
Fix $i$. 
Let $C_i\rhd D_i$ be the formula as given by Condition \ref{i:sdef:3}.
We thus have $C_i\rhd D_i\in(w)_1$, where $\Box\neg C_i\in Y_i$ and,
$C_i\in(x_i)_1$ or $\Diamond C_i\in(x_i)_1$.
For any $j\geq i$, this implies $\Box\neg C_i\in Y_j$ which gives
$\Box\neg C_i\in(y_j)_1$ and thus $\neg C_i,\Box\neg C_i\in(x_{j + 1})_1$.
The latter gives $C_i\neq C_{j + 1}$,
which is a contradiction since $\Phi$ is finite.
\end{proof}

We conclude the proof of Theorem \ref{thm:ilw-complete} by proving a truth lemma.

\begin{lemma}[Truth lemma]
For all $F\in\Phi$ and $w\in W$ we have $F\in(w)_1$ iff $w\Vdash F$.
\end{lemma}
\begin{proof}

We proceed by induction on $F$. 

The cases of the propositional variables and the connectives
are easily provable using properties of \mcs{}s and the $\Vdash$ relation.
So suppose $F = A\rhd B$.




$(\Rightarrow)$ Suppose we have $A\rhd B\in (w)_1$.
Then for all $v$ such that $wRv$ and $v\Vdash A$ we have to find a $u$ such that $v S_w u\Vdash B$
which, by the induction hypothesis, is equivalent to $B\in(u)_1$.
Consider such a $v$.
We have for some $S$ that $(v)_0=(w)_0\ast \langle S\rangle\ast\tau$ and 
$(w)_1\fsassuring{S}(v)_1$.
By the induction hypothesis we see that $A \in (v)_1$, so by Lemma \ref{lemm:wfin:defies} there is a \mcs $\Delta$ such that
$(w)_1\fsassuring{S\cup \{\Box\neg A\}}\Delta \ni B$. 
We take $u=\langle (w)_0\ast\langle S\cup\{\Box\neg A\} \rangle ,\Delta  \rangle$.
Now \ref{i:sdef:3} holds whence $v S_w u$.

$(\Leftarrow)$ 
Suppose that $A\rhd B\notin (w)_1$.
Then $\neg(A\rhd B)\in(w)_1$ whence by Lemma \ref{lemm:wfin:probs} there is a \mcs $\Delta$ such that
$(w)_1\fsassuring{\{\Box\neg A, \neg B\}}\Delta\ni A$.
Consider  $v'=\langle (w)_0\ast\langle\{\Box\neg A, \neg B\} \rangle,\Delta \rangle$.
We claim there is no $u'$ such that $v'S_w u'\Vdash B$.
Suppose otherwise. 
If $v'S_w u'$ because of $v'Ru'$ or $v' = u'$, then $(w)_1 \fsassuring{ \{ \neg B \} } (u')_1$ (possibly using the fact that $(w)_1 \fsassuring{ \{ \neg B \} } (v')_1 \fsassuring{} (u')_1$ implies $(w)_1  \fsassuring{ \{ \neg B \} } (u')_1$). 
Otherwise, both (3a) and (3b) hold.
Then $(u')_0 = (w)_0 \ast \langle T \rangle \ast \tau$ for some $T \supseteq \{ \Box \neg A, \neg B \}$. 
Thus $(w)_1  \fsassuring{ \{ \neg B \} } (u')_1$.
\end{proof}


\section{Labels and transitive closure}
\label{sec:seq}

The labelling that was considered in this paper so far was concerned with two or three worlds at a time. 
Due to the transitivity of $R$, labelling longer sequences often simplifies to labelling pairs or triples of worlds.
To give an example, suppose the ambient logic is \extil{M_0} and we have a sequence $\Sigma \prec_U \Gamma \prec_S \Delta \prec_T \Theta$.
In this sequence of four worlds there are two sequences of (immediately neighbouring) three worlds. The labelling lemma for \extil{M_0} tells us the following facts for these two sequences:
\[
    \Sigma \prec_{U \cup \Gamma_\emptyset^\Box} \Delta \prec_T \Theta
    \hskip 1em \text{and} \hskip 1em
    \Sigma \prec_U \Gamma \prec_{S \cup \Delta^\Box_\emptyset} \Theta,
\]
hence 
\[
    \Sigma \prec_{U \cup \Gamma_\emptyset^\Box \cup \Delta_\emptyset^\Box} \Theta
    \hskip 1em \text{and} \hskip 1em
    \Sigma \prec_{U \cup \Gamma_\emptyset^\Box} \Theta,
\]
and clearly the first fact is more informative.
Thus with \extil{M_0} labelling sequences of worlds with more than three worlds simplifies to labelling triples. 

In this section, we recall that labelling sequences in \illuka{R}-models indeed reduces to labelling triples of worlds.
The completeness of the logic \illuka{R} w.r.t.\ the ordinary Veltman semantics is still an open problem.
The fact that labels for this logic are compatible with transitive closures makes our labelling a good candidate for the step-by-step completeness proofs such as the \textit{construction method} \cite{jogo08}.

In the next section we deal with logics  whose labelling does not trivially reduce to labelling pairs or triples of worlds.
At the moment, the only logics falling into this category that we know of are various extensions of \illuka{W}. An example is \illuka{WR}, which may also be the most interesting example since it is the simplest logic among those whose (in)completeness status is currently open. Note that \extil{R} is not an extension of \extil{W}.
This is also the reason why we emphasise that for the logic \extil{R}, when taken on its own, the labelling is simple.

Suppose we are working in \illuka{R}. Let us recall the labelling lemma for \illuka{R}, Lemma \ref{lemm:labr}:
$\Gamma \sassuring{S} \Delta \sassuring{T} \Delta' \Rightarrow 
\Gamma \sassuring{S \cup \boxset{T}{\Delta}} \Delta'$.

Consider the situation as before with \extil{M_0}, but with the new ambient logic \extil{R}: $\Sigma \prec_U \Gamma \prec_S \Delta \prec_T \Theta$.
The labelling lemma for \extil{R} tells us the following facts for these two sequences:
\[
    \Sigma \prec_{U \cup \Gamma_S^\Box} \Delta \prec_T \Theta
    \hskip 1em \text{and} \hskip 1em
    \Sigma \prec_U \Gamma \prec_{S \cup \Delta^\Box_T} \Theta,
\]
hence 
\[
    \Sigma \prec_{U \cup \Gamma_S^\Box \cup \Delta_T^\Box} \Theta
    \hskip 1em \text{and} \hskip 1em
    \Sigma \prec_{U \cup \Gamma_{S \cup \Delta^\Box_T}^\Box} \Theta.
\]

Now it is not obvious which item is more informative,
    and the following lemma from \cite{Bilkova-Goris-Joosten} answers this question: the first fact is more informative.

\begin{lemma}\label{lemm:rtransstrat}
For logics containing 
 \principle{R} we have
$\Gamma \sassuring{S} \Delta  
 \Rightarrow
\boxset{S \cup \boxset{T}{\Delta}}{\Gamma} \subseteq \boxset{T}{\Delta}$.
\end{lemma}

\begin{proof}
See \cite{Bilkova-Goris-Joosten}.
\end{proof}

Thus, with \extil{R} and extensions, labelling sequences of worlds with more than three worlds simplifies to labelling triples.

\section{Non-trivial labellings of sequences}\label{sec:seq_gen}

We will first present an issue concerning labelling in \illuka{WR}. Both \illuka{W} \cite{jonvelt99} and \illuka{R} \cite{complgensem20} are known to be complete, but this question remains open for \illuka{WR}.
We will then proceed to work with \illuka{P}, another logic exhibiting the same issue (if we wish to prove a slightly stronger completeness result than the standard one). 
We switch from \illuka{WR} to \illuka{P} because we do not have a full proof of completeness of \illuka{WR} yet, while with \illuka{P} we can give a full completeness proof together with a to-the-point presentation on how to deal with logics with non-trivial labelling of sequences.

\tikzset{every picture/.style={line width=0.75pt}} 

\begin{figure}
    \centering
    
    \scalebox{0.8}{
        \begin{tikzpicture}[x=0.75pt,y=0.75pt,yscale=-1,xscale=1]
        
        \draw    (94.5,304.91) -- (215.5,304.91) ;
        \draw [shift={(217.5,304.91)}, rotate = 180] [color={rgb, 255:red, 0; green, 0; blue, 0 }  ][line width=0.75]    (10.93,-3.29) .. controls (6.95,-1.4) and (3.31,-0.3) .. (0,0) .. controls (3.31,0.3) and (6.95,1.4) .. (10.93,3.29)   ;
        \draw    (234.5,304.91) -- (342.41,304.91) ;
        \draw [shift={(344.41,304.91)}, rotate = 180] [color={rgb, 255:red, 0; green, 0; blue, 0 }  ][line width=0.75]    (10.93,-3.29) .. controls (6.95,-1.4) and (3.31,-0.3) .. (0,0) .. controls (3.31,0.3) and (6.95,1.4) .. (10.93,3.29)   ;
        \draw    (355.5,315) .. controls (371.26,340.61) and (334.63,337.11) .. (361.24,367.58) ;
        \draw [shift={(362.5,369)}, rotate = 227.82] [color={rgb, 255:red, 0; green, 0; blue, 0 }  ][line width=0.75]    (10.93,-3.29) .. controls (6.95,-1.4) and (3.31,-0.3) .. (0,0) .. controls (3.31,0.3) and (6.95,1.4) .. (10.93,3.29)   ;
        \draw    (370.46,382.09) -- (464.44,383.28) ;
        \draw [shift={(466.44,383.3)}, rotate = 180.72] [color={rgb, 255:red, 0; green, 0; blue, 0 }  ][line width=0.75]    (10.93,-3.29) .. controls (6.95,-1.4) and (3.31,-0.3) .. (0,0) .. controls (3.31,0.3) and (6.95,1.4) .. (10.93,3.29)   ;
        \draw    (94.5,304.91) -- (351.58,378.45) ;
        \draw [shift={(353.5,379)}, rotate = 195.96] [color={rgb, 255:red, 0; green, 0; blue, 0 }  ][line width=0.75]    (10.93,-3.29) .. controls (6.95,-1.4) and (3.31,-0.3) .. (0,0) .. controls (3.31,0.3) and (6.95,1.4) .. (10.93,3.29)   ;
        \draw  [dash pattern={on 4.5pt off 4.5pt}]  (369.1,304.91) .. controls (428.9,311.93) and (415.77,355.12) .. (464.93,382.48) ;
        \draw [shift={(466.44,383.3)}, rotate = 208.19] [color={rgb, 255:red, 0; green, 0; blue, 0 }  ][line width=0.75]    (10.93,-3.29) .. controls (6.95,-1.4) and (3.31,-0.3) .. (0,0) .. controls (3.31,0.3) and (6.95,1.4) .. (10.93,3.29)   ;
        \draw    (234.5,304.91) -- (464.55,382.66) ;
        \draw [shift={(466.44,383.3)}, rotate = 198.68] [color={rgb, 255:red, 0; green, 0; blue, 0 }  ][line width=0.75]    (10.93,-3.29) .. controls (6.95,-1.4) and (3.31,-0.3) .. (0,0) .. controls (3.31,0.3) and (6.95,1.4) .. (10.93,3.29)   ;
        \draw   (363,299.68) -- (421.5,225) -- (455.4,278.16) -- cycle ;
        
        \draw (82.2,303.29) node   [align=left] {$\displaystyle w$};
        \draw (225.09,302.48) node   [align=left] {$\displaystyle x$};
        \draw (355.89,304.48) node   [align=left] {$\displaystyle u$};
        \draw (362.89,378.67) node   [align=left] {$\displaystyle v$};
        \draw (473.18,381.67) node   [align=left] {$\displaystyle z$};
        \draw (152.36,294) node   [align=left] {$\displaystyle S$};
        \draw (294.07,295.61) node   [align=left] {$\displaystyle T$};
        \draw (208.73,352.02) node  [rotate=-16.99] [align=left] {$S \cup x^{\Box}_{T} \cup \{ \Box \neg C \}$ };
        \draw (342.49,353.87) node   [align=left] {$\displaystyle w$};
        \draw (453.29,357.39) node   [align=left] {$\displaystyle x$};
        \draw (294.89,333.04) node  [rotate=-20.77] [align=left] {$\displaystyle T$};
        \draw (410.07,266.61) node   [align=left] {$\displaystyle C$};
        \end{tikzpicture}
    }

    \caption{Labels with \illuka{WR}}
    \label{fig:ilwr_problem}
\end{figure}
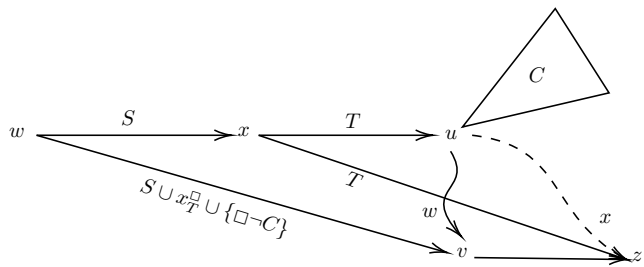

Recently, tools from an earlier version \cite{Bilkova-Goris-Joosten} of this paper have been used as one of the key ingredients in the completeness proof of \illuka{R} and other logics w.r.t.\ generalised semantics \cite{complgensem20}.
A natural next step is to tackle the completeness of \illuka{WR}.
However, when one tries to combine the completeness proofs for \illuka{W} and \illuka{R}, a problem occurs. 
Let us first discuss this problem, and then see how a more elaborate labelling system can help.
At the moment we do not know if the labelling systems will lead to a completeness proof of \illuka{WR}.
However the panorama looks promising. The problem that arises when addressing \illuka{WR} also occurs in a completeness proof for \illuka{P} and there we can solve it.
Thus, the above mentioned elaborate labelling systems should at least be an ingredient, if not the whole solution, in proving the more interesting completeness results.

Suppose\footnote{This paragraph describes the situation represented in Figure \ref{fig:ilwr_problem}.} we are building a model step-by-step (as in the \textit{construction method} \cite{jogo08}) and we have $A \rhd B \in w \prec_S x \prec_T u \ni A$.
So, we need to find some $v$ with $B \in v$ and a sufficiently strong label for $wRv$; and then declare $uS_w v$.
Using the labelling lemmas for \textsf{W} and \textsf{R}, it is easy to find $v$ with $w \prec_{S \cup x_T^\Box \cup \{ \Box \neg C \} } v$ for some $C$ contained either in $u$ or in a world $R$-accessible from $u$.
Let us for the moment suppose that any such $v$ fits our purposes.

Now, assume that at some later point during the construction, a world $z$ appears with $vRz$. 
By the frame condition of the principle \textsf{R}, we should have $uS_x z$.
If we were building an \illuka{R}-model (and not an \illuka{WR}-model), we would have to ensure just that $z$ has the same assuringness as $u$ with respect to $x$, that is, $xRz$ should be labelled with $T$.
Since we are building an \illuka{WR}-model and in order to ensure the frame condition for \textsf{W}, in addition to that we are to find a formula $C'$ with $x \prec_{T \cup \{ \Box \neg C' \} } z$.
An obvious candidate for $C'$ is $C$.
However, from $w \prec_{S \cup x_T^\Box \cup \{ \Box \neg C \} } v \prec z$ we only get $x \prec_{T} z$ (Lemma 22, \cite{complgensem20}), and what we would like is to have $x \prec_{T \cup \{ \Box \neg C \} } z$.
Let us refer to this phenomenon as the problem of label iteration.

One way to solve this problem is to simply \textit{require} $\Box \neg C$ to appear at the right place in the original label, i.e., instead of asking for $w \prec_{S \cup x_T^\Box \cup \{ \Box \neg C \} } v$, we ask for 
\begin{equation}
    \label{ilwp_labels_triple}
    w \prec_{S \cup x_{T \cup \{ \Box \neg C \}}^\Box \cup \{ \Box \neg C \} } v.
\end{equation}
If we are proving completeness w.r.t.\ generalised semantics using the approach from \cite{complgensem20}, this means that we should add a new condition in the definition of $S_w$ (Definition 28, \cite{complgensem20}). 
However, similar to how the original condition concerning two worlds requires us to add the new condition concerning three worlds that we just described, this condition itself requires us to add another condition, this time concerning four worlds. 
Let us illustrate this.

Suppose we have the following situation (see Figure \ref{fig:ilwr_problem}): \[
    A \rhd B \in w \prec_S x \prec_T u \ni A \text{ \  and \  } v \prec z.
\]
We would like to show $uS_x z$.
In particular, we have to show that if $x \prec_S u' \prec_T u$, then there is $v$ with 
    $x \prec_{ S \cup \{ \Box \neg A \} \cup (u')^\Box_{T \cup \{ \Box \neg A \}} } v \ni B$.
A good choice for such a world $v$ should satisfy \begin{equation}
    \label{eqn:wr_wanted}
    w \prec_{ x^\Box_{    S \cup \{ \Box \neg A \} \cup (u')^\Box_{T \cup \{ \Box \neg A \}}   } \cup \{ \Box \neg A \} } v,
\end{equation}
since from this we can conclude $x \prec_{ S \cup \{ \Box \neg A \} \cup (u')^\Box_{T \cup \{ \Box \neg A \}} } v \ni B$.
However, to be able to conclude \eqref{eqn:wr_wanted} we need to have a new case in the definition of $S_w$, one that concerns not just $u$ or just $x$ and $u$, but $x$, $u$ and $u'$.
Analogous reasoning applies for longer sequences of worlds.

It turns out the problem of label iteration, that, as we just saw, occurs with \illuka{WR}, also occurs when trying to prove that \illuka{P} is complete w.r.t.\ the class of generalised \illuka{P}-frames where an additional requirement which ensures \kgen{W} is present.\footnote{
    \kgen{W} is the following condition: $uS_w V \Longrightarrow (\exists V' \subseteq V) uS_w V' \mathbin{\&} R[V'] \cap S^{-1}_w[V] = \emptyset$.\\
    The requirement we mention is that whenever $w \prec_S u$ and we are making an $S_w$-successor $v$ of $u$, that $w \prec_{S \cup \{ \Box \neg B \}} v$ for some $B \in \mathcal{D} \cap \bigcup\dot{R}[u]$ where $\dot{R}[u] = R[u] \cup \{ u \} $. Since it is well known that $\illuka{P}$, which is a complete logic \cite{jovelt90}, contains $\illuka{W}$ (see e.g.\ \cite{Visser:1990:InterpretabilityLogic}), we already know that \illuka{P} is complete w.r.t.\ the class of generalised \illuka{P}-frames that satisfy \kgen{W}. 
    We do not, however, know in general if the models obtained by the standard completeness argument also satisfy this specific requirement (which is, at least a priori, stronger than \kgen W). 
}
In the next section we will give a detailed exposition on how to handle this problem in the case of \illuka{P}. 
The same general approach should be useful for any other extension of \illuka{W} that exhibits the problem of label iteration. 

\section{Completeness of \illuka{P} w.r.t.~a restricted class of frames}

In this section we will introduce the labelling system, i.e.~a systematic way to assign labels to arbitrarily long sequences, for \illuka{P}, and prove the completeness of \illuka{P} w.r.t.\ the class of generalised \illuka{P}-frames where an additional requirement which ensures \kgen{W} is present. This proof first appeared in the thesis \cite{luka-phd}.

\subsection{A labelling system}
The characteristic property of \kgen{P} is:
\[ wRw'RuS_w V \Rightarrow (\exists V' \subseteq V)\ uS_{w'} V'.  \]

Recall the labelling lemma for \illuka{P} (Lemma \ref{sec:labellinglemmata}.\ref{lemm:labp}):
\[
    w \prec_S x \prec_T u   \Rightarrow   w \prec_{S \cup x_T^\boxdot} u.
\]

The actual labelling that we use is an iterated generalisation of this property. 
Thus, instead of defining labels between pairs of MCSs, we consider tuples of MCSs with labels between them: $w_n \prec_{S_n} w_{n - 1} \prec_{S_{n - 1}} \dots \prec_{S_1} w_0 $.
We wish to define labels for \il{P} similar to the ones for \il{WR} between $w$ and $v$ in \eqref{ilwp_labels_triple} and \eqref{eqn:wr_wanted}.
We will first define these labels, and then prove the appropriate labelling lemma.

\begin{definition}
    \label{defn-Qn}
    For $n \in \omega {\setminus} \{0\}$, let $(w_i)_{ i \in \{ 0, \dots, n\}}$ be a finite sequence of \illuka{P}-MCSs, let $(S_i)_{ i \in \{ 1, \dots, n \}}$ be a finite sequence of sets of formulas and $B$ be a formula.
    We will define a sequence of sets of formulas, and we will denote this sequence as
    \[ 
        (Q((w_i)_{ i \in \{ 0, \dots, n\}}, (S_i)_{ i \in \{ 1, \dots, n \}}, B, j))_{ j \in \{ 1, \dots, n \}}.
    \] 
    Usually the MCSs $(w_i)_{ i \in \{ 0, \dots, n\}}$ and the sets of formulas $(S_i)_{ i \in \{ 1, \dots, n \}}$ will be clear from the context, so we will write $Q_j(B)$ for $Q((w_i)_{ i \in \{ 0, \dots, n\}}, (S_i)_{ i \in \{ 1, \dots, n \}}, B, j)$.
    We now recursively define the elements of our sequence:
    \begin{align*}
        Q_1(B)  &:= S_1 \cup \{ \Box \neg B \}; \\
        Q_{j + 1}(B) &:= S_{j + 1} \cup \{\Box \neg B \} \cup (w_{j})^\boxdot_{ Q_j(B) }.
    \end{align*}
\end{definition}

Note that the preceding definition amounts to the following: \[
    Q_j(B) = S_j \cup \{ \Box \neg B \} \cup (w_{j-1})^\boxdot_{ S_{j-1} \cup \{\Box \neg B \}  \cup (w_{j - 2})^\boxdot_{ {S_{j - 2} \cup \{\Box \neg B \}  \cup \dots}_{\dots \cup (w_1)^\boxdot_{ S_1 \cup \{\Box \neg B \} }  }}}.
\]

We now show that our labelling system is not too strong: we really can construct worlds satisfying such labels. In addition to this, we will later show that our labelling system is strong enough to ensure that the characteristic property holds in models that we build in the completeness proof. 

\begin{lemma}
\label{nedostaci-WR}
Let $n \in \omega {\setminus} \{0\}$ be arbitrary, $(w_i)_{ i \in \{ 0, \dots, n\}}$ be a finite sequence of \illuka{P}-MCSs, $(S_i)_{ i \in \{ 1, \dots, n \}}$ a finite sequence of sets of formulas and $B \rhd C$ a formula such that: \[
    B \rhd C \in w_{n} \prec_{S_{n }} w_{n - 1} \prec_{S_{n - 1}} \dots \prec_{S_{1}} w_0 \ni B.
\]
Then there is an \illuka{P}-MCS $v$ such that $w_n \prec_{Q_{n}(B)} v$ and $C, \Box \neg C\in v$.
\end{lemma}
\begin{proof}   
    We prove the claim by induction on $n$.
    In the base case we are to find $v$ such that $w_1 \prec_{S_1 \cup \{ \Box \neg B \} } v$. 
    But this is just Lemma \ref{lemm:wdefies}.
    
    Let us prove the claim for $n + 1$.
    Fix MCSs $(w_i)_{ i \in \{ 0, \dots, n, n + 1\}}$, labels $(S_i)_{ i \in \{ 1, \dots, n, n + 1\}}$ and a formula $B \rhd C$.
    Assume \[
        B \rhd C \in w_{n + 1} \prec_{S_{n + 1}} w_{n } \prec_{S_{n }} \dots \prec_{S_{1}} w_0\ni B.
    \]
    The goal is to find $v$ with $w_{n + 1} \prec_{Q_{n + 1}(B)} v \ni C, \Box \neg C$, i.e. \[
        w_{n + 1} \prec_{S_{n + 1} \cup \{ \Box \neg B \} \cup (w_n)^\boxdot_{ Q_{n }(B) }} v \ni C, \Box \neg C.
    \]
    From $w_{n + 1} \prec_{} w_{n}$ and the axiom \textsf P we have $B \rhd C \in w_{n }$.
    By the induction hypothesis, there is $v$ with $w_{n} \prec_{ Q_{n}(B) } v \ni C, \Box \neg C$.
    From $w_{n + 1} \prec_{S_{n + 1}} w_n \prec_{ Q_{n}(B) } v$ and the labelling lemma for \illuka{P} (Lemma \ref{lemm:labp}) we have: \[
        w_{n + 1} \prec_{S_{n + 1} \cup (w_n)^\boxdot_{ Q_{n }(B) }} v.
    \]
    Since $\{ \Box \neg B \} \subseteq Q_n(B) \subseteq (w_n)^\boxdot_{ Q_{n }(B) }$, we have $S_{n + 1} \cup (w_n)^\boxdot_{ Q_{n }(B) } = S_{n + 1} \cup \{ \Box \neg B \} \cup (w_n)^\boxdot_{ Q_{n }(B) }$.
\end{proof}
Note that the last line shows that a simpler definition of $Q_{j + 1}(B)$ would suffice: $Q_{j + 1}(B) := S_{j + 1} \cup (w_{j})^\boxdot_{ Q_j(B) }$ instead of $Q_{j + 1}(B) := S_{j + 1} \cup \{\Box \neg B \} \cup (w_{j})^\boxdot_{ Q_j(B) }$.
However, the purpose of this section is to introduce a method for dealing with arbitrary extensions of \illuka{W}.
We do not think it is likely that such a simplification could be made in the case of more interesting logics, such as \illuka{WR}.

\subsection{\extil{WP}-structures}

In the remainder of this section, $\mathcal D$ will always be assumed to be a finite set of formulas closed under taking subformulas and single negations, and $\top\in \mathcal{D}$ (i.e.\ $\bot \to \bot \in \mathcal{D}$).

Now we define the structures w.r.t.\ which we later prove completeness.
When defining $S_w$ we have to take care to make it compatible with the properties of a generalised Veltman model, in particular, the property that $wRu$ implies $uS_w \{ u \}$ and the property that $wRuRv$ implies $u S_w \{ v \}$.
So, if we fix $w$ and $u$, we should have $u S_w \{ v \}$ for all $v \in \dot{R}[u] (= R[u] \cup \{ u \})$.
However, because of monotonicity, we want not only $u S_w \{ v \}$ in such cases, but also $u S_w V$ for all $V \subseteq R[v]$ that contain $v$. This is why we add the condition (a) in the definition below.

Note that in the definition below, worlds are sets of formulas. Because of this, the operation $\bigcup\dot{R}[u]$ makes sense and defines a set of formulas.

\begin{definition}
We say that $\mathfrak{M} = (W, R, \{S_w : w \in W\}, \Vdash)$ is the \illuka{WP}-structure for a set of formulas
$\mathcal{D}$ if:
\begin{itemize}
    \item $W = \{ w : w \text{ is an \illuka{P}-MCS and for some } B \in \mathcal{D}, 
            \ B \wedge \square \neg B \in w \}$;
    \item $wRu \Leftrightarrow w \prec u$;
    \item $uS_w V \Leftrightarrow wRu \text{ and } V \subseteq R[w] \text{ and, moreover, one of the following 
            holds:}$
    
        $\displaystyle
            \begin{aligned} 
                (a) \ &V \cap \dot{R}[u] \neq \emptyset; \\
                (b) \ &\text{we have for all } 
                    n \in \omega {\setminus} \{0\}, 
                    \text{ all }
                    (w_i)_{ i \in \{ 0, \dots, n \}},
                    \text{ and all }
                    (S_i)_{ i \in \{ 1, \dots, n  \}}{:} \\
                & w = w_n \prec_{S_n} \dots \prec_{S_{1}} w_0 = u \Rightarrow (\exists v \in V) (\exists B \in 
                \mathcal{D} \cap \bigcup\dot{R}[u]) \ w \prec_{ Q_{n}(B) } v;
            \end{aligned}
        $ 
    \item $w\Vdash p \Leftrightarrow p\in w$.
\end{itemize}
\end{definition}


\begin{lemma} \label{lemma-main-wp}
The \illuka{WP}-structure $\mathfrak{M}$ for $\mathcal{D}$ is a generalised Veltman model. Furthermore, the 
following holds for each $w\in W$ and $G\in\mathcal{D}$:
$$ \mathfrak{M},w\Vdash G \ \mbox{ if and only if }\ G\in w,$$
\end{lemma}
\begin{proof}
Let us first verify that the \illuka{WP}-structure $\mathfrak{M} = (W, R, \{ S_w : w \in W \}, {\Vdash})$ 
for $\mathcal{D}$ is a generalised Veltman model. 
All the properties, except for quasi-transitivity, have easy proofs (see \cite{complgensem20}, the proof of Lemma 29).

Let us prove quasi-transitivity. Thus, we assume $u S_w V$, and $v S_w U_v$ for all $v \in V$. 
We put $U = \bigcup_{v \in V} U_v$ and claim that $u S_w U$. Clearly $U \subseteq R[w]$. 
To prove $uS_w U$ we will distinguish cases from the definition of the 
relation $S_w$ for $uS_w V.$

In Case (a), there exists a MCS $v_0 \in V$ for some $v_0 \in \dot{R}[u]$. 
We will next distinguish two Cases from the definition of $v_0 S_w U_{v_0}$. 

In Case (aa) we can find $x \in U_{v_0}$ for some $x \in \dot{R}[v_0]$. 
Since $v_0 \in \dot{R}[u]$, also $x \in \dot{R}[u]$. 
And since $x \in U_{v_0} \subseteq U$, we have $U\cap \dot{R}[u]\neq \emptyset.$
So, we have $uS_w U$ as required. 

In Case (ab):
\begin{align*}
    & \text{For all } n \in \omega {\setminus} \{0\}, \text{ all } 
    (w_i)_{ i \in \{ 0, \dots, n \}},
    \text{ and all }
    (S_i)_{ i \in \{ 1, \dots, n  \}}
    \text{ we have: }\\
    & w = w_n \prec_{S_n} \dots \prec_{S_{1}} w_0 = v_0 \Rightarrow (\exists x \in U_{v_0}) (\exists B \in 
    \mathcal{D} \cap \bigcup\dot{R}[v_0]) \ w \prec_{ Q_{n}(B) } x.
\end{align*}
To prove $uS_w U$ in this case, we will use Case (b) from the definition of the relation $S_w$. 
Let $n \in \omega {\setminus} \{0\}$ be arbitrary and let $(w_i)_{ i \in \{ 0, \dots, n \}}$ and $(S_i)_{ i \in \{ 1, \dots, n  \}}$ be arbitrary such that
    $w = w_n \prec_{S_n} \dots \prec_{S_{1}} w_0 = u$. 
If $u = v_0$, applying the formula above with the worlds $(w_i)_{ i \in \{ 0, \dots, n \}}$ and the labels $(S_i)_{ i \in \{ 1, \dots, n  \}}$ produces the required $x \in U_{v_0}$ and $B \in \mathcal{D} \cap \bigcup\dot{R}[v_0]$.
Otherwise, i.e.\ if $uRv_0$, let $w'_0 = v_0$, $w'_{i + 1} = w_i$, $S'_{1} = \emptyset$, $S'_{i + 1} = S_{i}$ and apply the formula above with $n + 1$, the sequence $(w'_i)_{ i \in \{ 0, \dots, n + 1 \}}$ and the labels $(S'_i)_{ i \in \{ 1, \dots, n + 1  \}}$.
This gives us a world $x \in U_{v_0}$ and a formula $B \in \mathcal{D} \cap \bigcup\dot{R}[v_0]$ with: \[
    w \prec_{ S_{n} \cup \{ \Box \neg B \} \cup (w_{n - 1})^\boxdot_{ {S_{n - 1} \cup \dots}_{\dots (w_1)^\boxdot_{ S_1 \cup \{ \Box \neg B \} \cup {u}^\boxdot_{ \emptyset \cup \{ \Box \neg B \} } }}}} x.
\]
Weakening this fact by Lemma \ref{lemm:trivialsassure} with removing ${u}^\boxdot_{ \emptyset \cup \{ \Box \neg B \} }$, we have the required property.
Since $uRv_0$ or $u=v_0$, we have $\dot{R}[v_0] \subseteq \dot{R}[u]$.
Thus we can reuse $B$ for this $S_w$ transition.

\vskip 2ex
In Case (b): 
\begin{align*}
    &\text{For all } n \in \omega {\setminus} \{0\}, \text{ all } 
    (w_i)_{ i \in \{ 0, \dots, n \}},
    \text{ and all }
    (S_i)_{ i \in \{ 1, \dots, n  \}}
    \text{ we have: } \\
    & w = w_n \prec_{S_n} \dots \prec_{S_{1}} w_0 = u \Rightarrow (\exists v \in V) (\exists B \in 
    \mathcal{D} \cap \bigcup\dot{R}[u]) \ w \prec_{ Q_{n}(B) } v.
\end{align*}

To prove $uS_w U$ we will use Case (b) from the definition of the relation $S_w$. 
So, let $n \in \omega {\setminus} \{0\}$ be arbitrary and let $(w_i)_{ i \in \{ 0, \dots, n \}}$ and $(S_i)_{ i \in \{ 1, \dots, n  \}}$ be arbitrary such that
    $w = w_n \prec_{S_n} \dots \prec_{S_{1}} w_0 = u$. 

By the assumption of this case, there are $v_0 \in V$ and $B \in \mathcal{D} \cap \bigcup\dot{R}[u]$ such 
that $w \prec_{Q_{n }(B)} v_0$. 
From $v_0 \in V$ we have $v_0 S_w U_{v_0}$. 
We will next distinguish the possible cases in the definition for $v_0 S_w U_{v_0}$. 

In the first Case (ba) we have $U_{v_0}\cap \dot{R}[v_0]\neq \emptyset,$ i.e.\ there is $x \in U_{v_0}$ such that either $v_0 = x$ or $v_0Rx$.
In both cases we have $w \prec_{ Q_n(B)} x$. 

In Case (bb), we have (Case (b) for $v_0 S_w U_{v_0}$ applied to $n = 1$ and $S_1 = Q_n(B)$) that
 there are some 
$x \in U_{v_0}$ and $B' \in \mathcal{D} \cap \bigcup\dot{R}[v_0]$ such that 
$w \prec_{Q_{n }(B) \cup \{ \square\neg B'\} } x$. 
By weakening (Lemma \ref{lemm:trivialsassure}), $w \prec_{Q_{n }(B)} x$, as required.

\vskip 3ex
We claim that for each formula $G\in\mathcal{D}$ and each world $w\in W$ the following holds:
\[
    \mathfrak{M},w\Vdash G \ \mbox{ if and only if } \ G\in w.
\]
The proof is by induction on the complexity of $G$. 
The only non-trivial case is when $G=B\rhd C.$

Assume $B \rhd C \in w,$ \ $wR u$ and $u\Vdash B$. Induction hypothesis implies $B\in u.$
We claim that $u S_w \{x : wRx\Vdash C \}$  by Case (b) from the definition of $S_w$. 
Clearly $wRu$ and $\{x : wRx\Vdash C \} \subseteq R[w]$. 
    
Fix $n \in \omega {\setminus} \{0\}$, $(w_i)_{ i \in \{ 0, \dots, n\}}$ and $(S_i)_{ i \in \{ 1, \dots, n  \}}$.
Assume $w = w_n \prec_{S_n} \dots \prec_{S_{1}} w_0 = u$. 
Since $B \rhd C \in w_n$ and $B \in w_0$, Lemma \ref{nedostaci-WR} implies that there is an \illuka{P}-MCS $v$ with $w_n \prec_{Q_{n}(B)} v$ and $C, \Box \neg C \in v$ (thus $v \in W$).
Since $C \in v$, the induction hypothesis implies $v \Vdash C$.
Since $w \prec v$, i.e.\ $wRv$, then $v \in \{x : wRx\Vdash C \}$. 
Finally, $B \in \mathcal{D}$ and $B\in u$ imply $B \in \mathcal{D} \cap \bigcup\dot{R}[u]$.

To prove the converse, assume $B \rhd C \notin w$. 
Since $w$ is an \illuka{P}-MCS,  $\neg(B\rhd C)\in w.$
Lemma \ref{lemm:wprobs} implies there is $u$ with 
$w \prec_{\{ \square \neg B, \neg C\}} u$ and  $B \in u.$ 
Since $w \prec_{\{ \square \neg B \}} u$, we have in particular that $\square \neg B\in u.$ 
So, $u \in W.$ 
The induction hypothesis implies $u \Vdash B$.
Let $V\subseteq R[w]$ be such that $u S_w V$.
We will find a world $v\in V$ such that $w\prec_{\{ \neg C\}} v$.  
We will distinguish Cases (a) and (b) from the definition of the relation $S_w$. 
Consider Case (a). 
Let $v$ be an arbitrary world in $V\cap \dot{R}[u]$. If $v = u$, clearly 
$w\prec_{\{\square\neg B,\neg C\}} v$.
If $uRv$, then we have $w\prec_{\{ \square\neg B,\neg C\}} u\prec v.$ 
This implies $w\prec_{\{ \square\neg B,\neg C\}} v$. 
Consider Case (b). 
From $w\prec_{\{ \square\neg B,\neg C\}} u$ and the definition of $S_w$ it follows that there is $v \in V$ such that (for some formula $D$) $w \prec_{\{\square \neg B, \neg C,\square\neg D\}} v$.
In both cases we have $w\prec_{\{ \neg C\}} v$; thus $C \notin v$.
Induction hypothesis implies $v \nVdash C$; whence $V \nVdash C$, as required.
\end{proof}

\begin{theorem}
\illuka{P} is complete w.r.t.\ the class of all generalised Veltman frames satisfying \kgen P. In particular, \illuka{P} is complete w.r.t.\ the class of \illuka{WP}-structures generated by all appropriate sets $\mathcal D$ (finite sets containing $\top$ that are closed under subformulas and single negations).
\label{tm-ilwr}
\end{theorem}
\begin{proof}
In the light of Lemma \ref{lemma-main-wp}, it suffices to show that the \illuka{WP}-structure $\mathfrak{M}$ for $\mathcal{D}$ possesses the property \kgen{P}.\footnote{
    In general, for example with the logic \illuka{WR}, we would want to verify if \kgen{W} holds.
    The proof would be the same as the proof of Theorem 30 in \cite{complgensem20}.
    In this case it is a consequence of \kgen{P}.
}

Let us prove \kgen{P}. Let $wRw'RuS_w V$ and take $V' = V \cap R[w']$. We claim $uS_{w'} V'.$ 

We distinguish two possible cases for $uS_w V$. 
If it holds by Case (a), there is $v \in V$ such that either $u = v$ or $uRv$.
In both cases $w'Rv$. 
Let $U = \{v\}$.
Clearly $U \subseteq V$.
Since $w'RuRv$, $uS_{w'} \{v\}$, i.e.\ $uS_{w'} U$.
The remainder of the proof deals with the case when $uS_w V$ holds by Case (b) from the definition of $S_w$.

Fix $n \in \omega {\setminus} \{0\}$, the worlds $(w_i)_{ i \in \{ 0, \dots, n\}}$ and the labels $(S_i)_{ i \in \{ 1, \dots, n  \}}$.
Assume $w' = w_n \prec_{S_n} \dots \prec_{S_{1}} w_0 = u$. 
We have $w \prec_\emptyset w_n \prec_{S_n} \dots \prec_{S_{1}} w_0$. 
Now the definition of $u S_w V$ implies there is $v \in V$ with:
\[
    w \prec_{ \emptyset \cup \{ \Box \neg B \} \cup (w_n)^\boxdot_{ Q_n (B)}} v.
\]
We claim that $w_n \prec_{Q_n(B)} v$.
Assume $\neg A \rhd \bigvee \neg F_i \in w_n$ with $F_i \in Q_n(B)$ (we are to show that $A, \Box A \in v$). Clearly $A, \Box A \in (w_n)^\boxdot_{ Q_n(B) }$. 
Since $a \prec_S b$ implies $S \subseteq b$, we have $A, \Box A \in v$.
\end{proof}

\bibliographystyle{alpha}
\bibliography{myref}

\end{document}